\documentclass[10pt]{article}
\usepackage{graphicx}
\usepackage{amsmath,amssymb,amsthm,amsfonts}
\usepackage{color}
\usepackage{cite}
\usepackage{amssymb}

\newtheorem{thm}{Theorem}[section]

\newtheorem{lem}[thm]{Lemma}
\newtheorem{prop}[thm]{Proposition}

\theoremstyle{definition}
\newtheorem{defn}{Definition}[section]
\newtheorem{ex}{Example}[section]
\theoremstyle{remark}
\usepackage{appendix}

\numberwithin{equation}{section}

\DeclareMathSymbol{\C}{\mathalpha}{AMSb}{"43}

\newcommand{\eps}{\varepsilon}

\newcommand{\lam}{\lambda}
\newcommand{\alp}{\alpha}
\newcommand{\R}{{\mathbb{R}}}
\newcommand{\h}{{\mathcal{H}}}
\newcommand{\inte}{\int_{\mathbb{R}^2}}

\newcommand{\f}{\alpha_a}

\def\o{\Omega}

\def\R{{\mathbb R}}
\def\C{{\mathbb C}}

\newcommand{\bsub}{\begin{subequations}}
\newcommand{\esub}{\end{subequations}$\!$}

\begin{document}	
\title{Local Uniqueness of Ground States for Rotating Bose-Einstein Condensates with Attractive Interactions}
\author{Yujin Guo\thanks{School of Mathematics and Statistics, and Hubei Key Laboratory of
Mathematical Sciences, Central China Normal University, P.O. Box 71010, Wuhan 430079, P. R. China. Email: \texttt{yguo@mail.ccnu.edu.cn}. Y. Guo is partially supported by NSFC under Grant  No. 11931012.
	},
	\, Yong Luo\thanks{School of Mathematics and Statistics, and Hubei Key Laboratory of
Mathematical Sciences, Central China Normal University, P.O. Box 71010, Wuhan 430079,
		P. R. China.  Email: \texttt{yluo@mail.ccnu.edu.cn}. Y. Luo is partially supported by the Project funded by China Postdoctoral Science Foundation No. 2019M662680.}
	\, and\, Shuangjie Peng\thanks{School of Mathematics and Statistics, and Hubei Key Laboratory of
Mathematical Sciences, Central China Normal University, P.O. Box 71010, Wuhan 430079,
		P. R. China.  Email: \texttt{sjpeng@mail.ccnu.edu.cn}. S. Peng is partially supported by the Key Project of NSFC under Grant No.11831009.
	}
}

\date{\today}

	\smallbreak \maketitle

\begin{abstract}
We study ground states of two-dimensional Bose-Einstein condensates  with attractive interactions in a trap  $V(x)$ rotating at the velocity $\Omega $.
It is known that there exists a critical rotational velocity $0<\Omega ^*:=\Omega^*(V)\leq \infty$ and a critical number $0<a^*<\infty$ {such that} for any rotational velocity $0\le \Omega <\Omega ^*$, ground states exist if and only if the coupling constant $a$ satisfies $a<a^*$. For a general class of traps $V(x)$, which may not be  symmetric, we prove in this paper  that up to a constant phase, there exists a unique ground state as $a\nearrow a^*$, where $\Omega\in(0,\Omega^*)$ is fixed.
\end{abstract}	
	
\vskip 0.05truein


\noindent {\it Keywords:} Bose-Einstein condensates; rotational velocity; local uniqueness; Pohozaev identity

\vskip 0.2truein

\section{Introduction}
A Bose-Einstein condensate (BEC) is a state of matter, in which atoms or particles are cooled to  the sufficiently low temperature that a large fraction of them ``condense" into a single quantum state. Because a BEC can present quantum effects at the macroscopic scale, it has become an important subject in experimental investigations since the first realization \cite{Anderson,D} of BECs in dilute gases of alkali atoms in 1995. Various interesting quantum phenomena have been observed in the physical experiments of BECs over the past two decades, including the critical-mass collapse \cite{Hulet1,Hulet2,HM,KM,Hulet3,D}, the appearance of quantized vortices \cite{Abo,B,CC,F,WG,ZW},  the center-of-mass rotation \cite{Abo,LC,F,WG}, and so on.
These novel experimental progresses promote greatly the developments of mathematical theories and numerical methods arising from BECs, see \cite{Abo,A,CR,D,F,IM-1,IM-2,Lieb06,LSS,LSY,WG}.

When the interactions between cold atoms in the condensates are repulsive, the quantized vortices and some other complex structures of BECs in rotating traps were analyzed and simulated extensively in the past few years, see review papers \cite{A2,LSS,Roug} and the references therein. Due to the distinct mechanism of the critical-mass collapse, the systems of rotating BECs in the attractive case however behave quite different from those of the well-understood repulsive case.
For example, vortices are generally unstable in rotating BECs with attractive interactions (see, e.g., \cite{CC,WG,LC}),
even though vortices are known cf.  \cite{A,CO,F} to form stable lattice configurations in the repulsive case.


By a mean-field approximation, see \cite{Lewin2,Lewin,LSS,Nam,Roug},  the energy of the two-dimensional attractive BEC in a rotating trap  can be described by the following Gross-Pitaevskii (GP) energy functional:
\begin{equation}
F_a(u):=\int _{\R ^2} \big(|\nabla u|^2+V(x)|u|^2\big)dx-\frac{a}{2}\int _{\R ^2}|u|^4dx-
\Omega \int_{\R ^2}x^{\perp}\cdot (iu,\, \nabla u)dx,  \ \ u\in \h , \label{f}
\end{equation}
where $x^{\perp} =(-x_2,x_1)$ with $x=(x_1,x_2)\in \R^2$,
$(iu,\, \nabla u)=i(u\nabla \bar u-\bar u\nabla u)/2$, and the complex space $\h$ is defined as
\begin{equation}\label{1.H}
\h :=  \Big \{u\in  H^1(\R ^2, \mathbb{C}):\ \int _{\R ^2}
V(x)|u|^2 dx<\infty\Big \}.
\end{equation}
Here the parameter $a>0$ in (\ref{f}) characterizes the absolute product of the scattering length $\nu$ of the two-body
interaction times the number $N$ of particles in the condensates, and $\Omega \ge 0$ describes the rotational velocity of the rotating trap $V(x)\ge 0$. Following (\ref{f}), {\bf ground states} of two-dimensional attractive BEC in a rotating trap satisfy (cf.\cite{BC,Lewin}) the following mass constraint  variational problem:
\begin{equation}\label{def:ea}
	e_F(a):=\inf _{\{u\in \h, \, \|u\|^2_2=1 \} } F_a(u),\ \, a>0.
\end{equation}
Alternatively, one may impose $e_F(a)$ a different constraint $\int_{\mathbb{R}^2} |u(x)|^2dx=N>0$, but the latter case can be easily  reduced to the previous one with $a$ being replaced by $a/N$. In view of this fact, in this paper we focus on the form of $e_F(a)$ instead. We remark that $e_F(a)$ is essentially a mass-critical constraint variational problem. The mass-subcritical case of $e_F(a)$, where the nonlinear term $|u|^4$ is replaced by $|u|^p$ for $2<p<4$, was studied as early as in the pioneering work of Esteban-Lions \cite{EL}.

For the non-rotational case $\Omega=0$ of $e_F(a)$,  the existence, stability and mass concentration of real-valued  minimizers were studied recently in \cite{GLW,GS,GWZZ,GZZ,Z} and the references therein.
In this case, it was shown in \cite{GLW,GS,GWZZ,GZZ} that
$e_F(a)$ admits real-valued minimizers if and only if $a<a^*$, where $a^*=\|w\|^2_{L^2(\R^2)}$ and $w=w(|x|)>0$ is the unique (cf. \cite{K,W}) positive radial solution of the following nonlinear scalar field equation
\begin{equation}
-\Delta u+ u-u^3=0\  \mbox{  in } \  \R^2,\,\ u\in H^1(\R ^2,\R).  \label{Kwong}
\end{equation}
By the analytical approach of \cite[Theorem II.1]{CL}, this further implies the following existence and nonexistence: for the non-rotational case $\Omega=0$, $e_F(a)$ admits
complex-valued minimizers (i.e., ground states) if and only if $a<a^*$.
Moreover, the mass concentration, symmetry breaking and other analytical properties of minimizers for $e_F(a)$ at $\Omega =0$ were analyzed recently for different types of $V(x)$, see \cite{GLW,GS,GWZZ,GZZ,LPY} and the references therein.


The rotational case $\Omega>0$ of $e_F(a)$  was discussed more recently in \cite{BC,ANS,GLY,Lewin}, where the existence, stability and the limit behavior of complex-valued minimizers were studied. Specially, if the general trapping potential $0\leq V (x)\in L^{\infty}_{loc}(\R^2)$ satisfies
\begin{equation}\label{A:V}
\underline{\lim} _{|x|\to\infty }\frac{V(x)}{|x|^2}>0,
\end{equation}
the following critical rotational velocity $\Omega ^*:=\Omega ^*(V)$ is defined in \cite{GLY}:
\begin{equation}
\Omega ^*:=\sup \Big\{\Omega >0:\ \  V(x)-\frac{\Omega ^2}{4}|x|^2  \to\infty \,\ \mbox{as}\,\
|x|\to\infty \Big\}.  \label{Omega}
\end{equation}
One can note that depending on $V(x)$, both $0<\Omega ^*<\infty$ and $\Omega ^*=\infty$ can happen. Under the assumption (\ref{A:V}), the following existence and non-existence of minimizers for $e_F(a)$ were proved in \cite[Theorem 1.1]{GLY}:

\vskip 0.05truein
\noindent{\bf Theorem A.} {\em Assume $V(x)\in L^\infty_{\rm loc}(\R^2)$ satisfies (\ref{A:V}) such that $\Omega^*\in (0,+\infty]$ in \eqref{Omega} exists. Then we have
\begin{enumerate}
		\item If $0\le\Omega <\Omega ^*$ and $0\leq a< a^*:=\|w\|^2_2$, then there exists at least
		one minimizer of $e_F(a)$.
		\item If $0\le\Omega <\Omega ^*$ and $a \ge a^*:=\|w\|^2_2$, then there
		is no minimizer of $e_F(a)$.
		\item If $\Omega >\Omega ^*$, then for any $a\ge 0$, there is no minimizer of $e_F(a)$.
\end{enumerate}}

If $e_F(a)$ admits a minimizer $u_a$, then the variational theory yields that  $u_a$ is a ground state of
the following Euler-Lagrange  equation
\begin{equation}
-\Delta u_a+V(x)u_a+i\, \Omega \, (x^{\perp}\cdot \nabla u_a)=\mu u_a+a|u_a|^2u_a\quad \mbox{in}
\, \  \R^2,  \label{eqn}
\end{equation}
where $\mu = \mu (a, \Omega, u_a )\in \R$ is a suitable Lagrange multiplier. We remark that  the complex-valued solutions of \eqref{eqn} were studied directly in \cite{Cao,EL,Peng,LPW} and the references therein, where the existence, uniqueness, and other analytical properties of complex-valued solutions were obtained via the Lyapunov-Schmidt reduction, topological degree method, and some other arguments. As for $e_F(a)$, suppose $V(x)$ satisfies the assumption (\ref{A:V}) such that $\Omega^*\in (0,+\infty]$ exists, and let $0\le\Omega <\Omega ^*$ be fixed, so that $V_\Omega (x):=V(x)-\frac{\Omega^2}{4}|x|^2\ge 0$ in $\R^2$. By deriving energy estimates and applying the elliptic PDE theory, under above assumptions it was proved in \cite{GLY,Lewin} that the minimizer $u_a$ of $e_F(a)$ concentrates at a global minimum point of $V_\Omega (x)$ as $a\nearrow a^*$, in the sense that  for some $\eps_a \searrow 0$ as $a\nearrow a^*$,
\begin{equation}
\eps_a u_a\big(\eps_a(x-x_0)\big)\to \frac{w(x)}{\sqrt {a^*}}\ \  \mbox{strongly in $H^1(\R^2,\,\mathbb{C})$ as}\,\ a\nearrow a^*,
\label{intro:lower}
\end{equation}
where $x_0\in\R^2$ is a global minimum point of $V_\Omega (x)$, i.e., $V_\Omega (x_0): =\inf _{x\in \R^2}V_\Omega (x)$.

Whether a physics system admits a unique ground state or not is a fundamental and interesting problem.
Based on the convergence (\ref{intro:lower}), it was further obtained in \cite[Theorem 1.3]{GLY} that  the minimizers of $e_F(a)$ must be unique and vortex-free  as $a\nearrow a^*$ for the harmonic case where $V(x)=|x|^2$. It deserves emphasis that the conclusions of \cite[Theorem 1.3]{GLY} were proved by the so-called method of inductive symmetry, which cannot however be extended to the non-radially symmetric case of $V(x)$. On the other hand, we should mention that the non-radially symmetric trap $ V(x)=|x|^2+|x_1|^2+\lam |x_2|^2$, where $x=(x_1, x_2)\in\R^2$ and $\lam >0,$  was already used in BEC experiments, see \cite{MC1,MC2} and the references therein. It is therefore natural to wonder whether the uniqueness of minimizers for $e_F(a)$ holds for the case where the trap $V(x)$ is not radially symmetric and however satisfies some additional assumptions.

Stimulated by above facts, the main purpose of this paper is to address the uniqueness of minimizers for $e_F(a)$ under a more general class of traps $V(x)$, which may not be symmetric. For this purpose, we now introduce the following homogeneous functions:

\begin{defn}
A function $h(x):\R^2\longmapsto \R$ is called homogeneous of degree $p\in\R^+$ (about the origin), if
\begin{equation}\label{1:V}
h(tx)=t^ph(x)\,\ \hbox{for any $t\in\R^+$ and $x\in \R^2$.}
\end{equation}
\end{defn}
\noindent One can note that if $\lim_{|x|\to\infty} h(x) = +\infty$, then $x=0$ is the unique minimum point of $h(x)$.
Following the above definition, we next assume that $V(x)$ satisfies
\begin{enumerate}
\item [\rm($V$).]
$V_{\Omega}(x):=V(x)-\frac{\Omega^2}{4}|x|^2\geq0$ and $\{x\in\R^2:\,V_{\Omega}(x)=0\}=\{0\}$, where $V_{\Omega}(x)$ satisfies
\begin{equation}\label{s.3}
|V_{\Omega}(x)|\le C e^{\gamma|x|}\quad\hbox{and}\ \ |\nabla V_{\Omega}(x)|\le C e^{\gamma|x|}\quad\hbox{for some $\gamma>0$ as $|x|\to\infty$},
\end{equation}
and
\begin{equation}\label{s.4}
V_{\Omega}(x)=h(x)+o(|x|^p)\ \, \hbox{and}\,\ \frac{\partial V_{\Omega}(x)}{\partial{x_j}}=
\frac{\partial h(x)}{\partial{x_j}}+o(|x|^{p-1})\,\ \hbox{as}\,\  |x|\to 0,\,\ j=1,2
\end{equation}
for some homogeneous function $0\le h(x)\in C^{1}(\R^2)$ of degree $1<p\le 2$, where $h(x)$ satisfies $\lim_{|x|\to \infty}h(x)=+\infty$.
\end{enumerate}



\begin{ex}\label{1:example}
Consider the following non-radially symmetric potential
\begin{equation}\label{t-V}
V(x)=|x|^2+a_1|x_1|^p+a_2|x_2|^p,\quad \hbox{$a_1\ge 0,\,  a_2\ge 0$ and $1<p\leq 2$,}
\end{equation}
where $x=(x_1, x_2)\in\R^2.$ We remark that when  $p=2$, (\ref{t-V}) gives the harmonic potential considered in \cite{MC1,MC2} for BEC experiments. One can  easily check that $V(x)\in C^{1,p-1}_{loc}(\R^2)$ satisfies the assumptions (\ref{A:V}) and $(V)$ for the homogeneous function $0\le h(x)\in C^{1}(\R^2)$ satisfying
\begin{equation}\label{s.7}
\left\{\begin{array}{lll}
h(x)=a_1|x_1|^p+a_2|x_2|^p,\,\ &\hbox{if}\,\  1<p<2;\\[3mm]
h(x)=\big(1-\frac{\Omega^2}{4}\big)|x|^2+a_1|x_1|^2+a_2|x_2|^2,\,\ &\hbox{if}\,\  p=2,
\end{array}\right.
\end{equation}
where $a_1\ge 0$ and $a_2\ge 0$ are as in (\ref{t-V}).
Moreover, one can check from (\ref{s.7}) that
\begin{equation}\label{deH}
H(y):=\inte h(x+y)w^2(x)dx
\end{equation}
admits a unique critical point $y_0=0$, which is also non-degenerate in the sense that
\begin{equation}\label{s.8}
\hbox{det}\Big(\displaystyle\inte \frac{\partial h(x)}{\partial x_j}\frac{\partial w^2(x)}{\partial x_l}dx\Big)_{j,l=1,2}\not =0.
\end{equation}
\end{ex}

Under above assumptions, the main result of this paper is concerned with the following local uniqueness.

\begin{thm}\label{thm1.2}
Suppose $V(x)\in C^{1,\alp}_{loc}(\R^2)$ $(0<\alp<1)$  satisfies (\ref{A:V}) and $(V)$ for some homogeneous function $0\le h(x)\in C^{1}(\R^2)$ of degree $p\in (1,2]$, and let $\Omega\in (0,\Omega^*)$ be fixed, where $\Omega^*>0$ is defined as in (\ref{Omega}).
Assume that $H(y)$ defined by \eqref{deH} has a unique non-degenerate critical point $y_0$,
then up to a constant phase, there exists a unique complex-valued minimizer of $e_F(a)$ for $a^*-a>0$ small enough.	
\end{thm}


We remark that the local uniqueness, up to a constant phase,  of Theorem \ref{thm1.2} holds in the following sense: there exists a minimizer  $u_a$ of $e_F(a)$ such that any minimizer $U_a$ of $e_F(a)$ satisfies $U_a\equiv u_ae^{i\theta_a}$ in $\R^2$ for $a^*-a>0$ small enough, where $\theta _a\in[0,2\pi)$ is a suitable constant phase depending on $a$. We note that depending on the shape of $h(x)$, the unique non-degenerate critical point of $H(y)$ can be nonzero.   Example \ref{1:example} gives an explicit example that the uniqueness of Theorem \ref{thm1.2} holds for non-radially symmetric traps $V(x)$.



Even though the similar local uniqueness of complex-valued solutions was investigated in \cite{Cao,GLY} and somewhere else, to our best knowledge, those arguments do not work in our situation. Actually, it seems that the topological degree argument (e.g. \cite{Cao}) does not work for proving Theorem \ref{thm1.2}, due to the possible multiplicity of the Lagrange multiplier $\mu = \mu (a)\in \R$ in the following Euler-Lagrange equation
\begin{equation}
-\Delta u_a+V(x)u_a+i\, \Omega \, (x^{\perp}\cdot \nabla u_a)=\mu u_a+a|u_a|^2u_a\ \ \mbox{in}
\, \  \R^2  \label{A:eqn}
\end{equation}
for any given $0<a<a^*$.
On the other hand, the local uniqueness of \cite[Theorem 1.3]{GLY} follows strongly from the conclusion that if $V(x)=|x|^2$, then $Im(u_a)\equiv 0$ as  $a\nearrow a^*$, which is however false generally for the non-radially symmetric case of $V(x)$.
Therefore, it is necessary to investigate a different approach for proving Theorem \ref{thm1.2}.
Motivated by \cite{CLL,Grossi,Deng,GL,GLW} and the references therein, we shall prove Theorem \ref{thm1.2} by constructing various Pohozaev identities, which were widely  used  in the existing literature of studying the real-valued elliptic PDEs.


To illustrate our main idea of proving Theorem \ref{thm1.2}, by contradiction we now suppose that up to a constant phase, there exist two different minimizers $u_{1,a}$ and $ u_{2,a}$ of $e_F(a)$ as $a\nearrow a^*$, in the sense that $u_{1,a}\not\equiv u_{2,a}e^{i\theta}$ for any constant phase $\theta =\theta(a)\in[0,2\pi)$. Motivated by \cite{CLL,Grossi,Deng,GL,GLW,GLY}, we make the following transformation  of $u_{j,a}$:
\begin{equation}\label{AA:D}
\tilde{u}_{j,a}(x):=\alp_a u_{j,a} \big(\alp_a(x+y_0)\big)e^{-i(\frac{\alp_a^2\Omega }{2} x\cdot y_0^{\perp})}e^{i\varphi_{j,a}}=  R_{j,a}(x)+i  I_{j,a}(x),\ \ j=1,2,
\end{equation}
where $\alp_a:=\frac{(a^*-a)^{\frac{1}{2+p}}}{\lambda}>0$ is given in (\ref{2.1}), the point $y_0\in\R^2$ is as in Theorem \ref{thm1.2},  and the constant phase $\varphi_{j,a}\in [0,2\pi )$
can be chosen properly such that
\begin{equation}\label{AA:E}
\inte w(x) I_{j,a}(x)dx=0,\ \ j=1,2.
\end{equation}
Under the assumptions of Theorem \ref{thm1.2}, we shall prove in Proposition \ref{5:prop} that  $\tilde u_{j,a}(x)$ satisfies
\begin{equation}
\lim_{a\nearrow a^*}\tilde u_{j,a}(x)=\frac{w(x) }{\sqrt{a^*}}\ \ \hbox{strongly\,\ in\,\ $L^\infty (\R^2, \mathbb{C})\cap H^1 (\R^2, \mathbb{C})$},\ \ j=1,2,
\label{N1.69}
\end{equation}
which is the first step of proving Theorem \ref{thm1.2}.

Because $u_{1,a}\not\equiv u_{2,a}e^{i\theta}$ for any constant phase $\theta =\theta(a)\in[0,2\pi)$, we next define
\begin{equation}\label{AA:F}
\eta_a(x):=\frac{\tilde{u}_{2,a}(x)-\tilde{u}_{1,a}(x)}{\|\tilde{u}_{2,a}-\tilde{u}_{1,a}\|_{L^\infty(\R^2)}}= \eta_{1,a}(x)+i \eta_{2,a}(x),
\end{equation}
where $ \eta_{1,a}$ and $\eta_{2,a}$ denote the real and imaginary parts of $\eta_a$, respectively.
In order to continue the proof of Theorem \ref{thm1.2}, the second step is to establish the refined $L^\infty$-uniform estimates of both $ \eta_{1,a}$ and $\eta_{2,a}$ as $a\nearrow a^*$. By making full use of (\ref{AA:E}), we shall reach this aim by proving Lemma \ref{1.23} and Proposition \ref{p-3.3}, which reveal  that $ \eta_{1,a}$ is the dominant  part of $\eta_a(x)$  as $a\nearrow a^*$. As a result, we shall be able to prove that up to a subsequence  if necessary, $(\eta_{1,a},\eta_{2,a})\to (\eta_1,\eta_2)$ uniformly in $L^\infty(\R^2)$ as $a\nearrow a^*$,  where $(\eta_{1},\eta_{2})$ satisfies the following system
\begin{equation}\label{new1}
\left\{
\begin{aligned}
-\Delta \eta_1+\eta_1-3w^2\eta_1&=-\Big(\frac{2}{a^*}\int_{\R^2}w^3\eta_1dx\Big)w
\ \ \mbox{in}\,\ \R^2,\\[2mm]
-\Delta \eta_2+\eta_2-w^2\eta_2&=0
\,\ \mbox{in}\,\ \R^2,\,\ \inte w\eta_2dx=0,\\
\end{aligned}
\right.
\end{equation}
and hence
\begin{equation}\label{AA:30}
\eta_1=b_0(w+x\cdot\nabla w)+\sum_{i=1}^2b_i\frac{\partial w}{\partial x_i}\quad\hbox{and}\quad \eta_2\equiv 0\ \ \hbox{in}\,\ \R^2
\end{equation}
for some constants $b_0$, $b_1$ and $b_2$.

In the third step of proving Theorem \ref{thm1.2}, we shall prove that $b_0=b_1=b_2=0$ holds for (\ref{AA:30}) by constructing various Pohozaev identities.
Due to the appearance of the rotating term $i(x^{\perp}\cdot\nabla \tilde u_{j,a})$, it however seems difficult to establish directly Pohozaev identities of the complex-valued function $\tilde u_{j,a}$. To overcome this difficulty, as illustrated before (\ref{3.34}), we shall construct various Pohozaev identities for the real part $R_{j,a}$ of $\tilde u_{j,a}:=  R_{j,a} +i  I_{j,a}$.  A key point of this strategy is to derive the refined estimates for the terms produced by the rotation. Applying Lemma \ref{prop-1} and Proposition \ref{p-3.3}, we shall prove that those terms produced by the rotation are lower order as $a\nearrow a^*$, based on which we shall be able to establish the system (\ref{3.34}) concerning the relationship between the homogeneous potential $h(x)$ and the coefficients $b_0$, $b_1$ and $b_2$. By deriving another type of Pohozaev identities, we shall further prove that $b_0=0$.
Following this fact and the non-degeneracy assumption of $y_0$, we shall derive from (\ref{3.34}) that $ b_1=b_2=0$ holds for (\ref{AA:30}), and hence,
\begin{equation}\label{K:30}
\eta _a=\eta_{1,a}+i\eta_{2,a} \to  \eta_0=\eta_1+i\eta_2\equiv 0 \ \ \mbox{uniformly in}\ \ L^\infty(\R^2)\ \ \mbox{as} \ \ a\nearrow a^*.
\end{equation}
On the other hand,  one can conclude from Lemma \ref{1.23} that $ \eta_0=\eta_1+i\eta_2\not\equiv 0$ in view of the fact that $\|\eta _a\|_{L^\infty}\equiv 1$. This is a contradiction, and Theorem \ref{thm1.2} therefore follows.

This paper is organized as follows. In Section 2, we shall analyze the $L^\infty$-uniform estimates of minimizers for $e_F(a)$ as $a\nearrow a^*$. In Section 3, we shall first derive a crucial gradient estimate of \eqref{1.26}, based on which  the refined limit profiles  of minimizers as $a\nearrow a^*$ are then established in Lemma \ref{1.23} and Proposition \ref{p-3.3}. Following the estimates of Section 3, the complete proof of Theorem \ref{thm1.2} is finally addressed in Section 4 by constructing various  Pohozaev identities.



\section{$L^\infty$-uniform estimates as $a\nearrow a^*$}
The purpose of this section is to address $L^\infty-$uniform estimates of the complex-valued minimizers for $e_F(a)$ as $a\nearrow a^*$.
Towards this purpose, we first introduce the following
Gagliardo-Nirenberg inequality
\begin{equation}\label{GNineq}
\inte |u(x)|^4 dx\le \frac 2 {a^*} \inte |\nabla u(x) |^2dx \inte |u(x)|^2dx ,\
\  u \in H^1(\R ^2, \R),
\end{equation}
where the equality is attained  (cf. \cite{W}) at the unique positive radial solution $w$ of \eqref{Kwong}. Moreover, it follows from \cite[Lemma 8.1.2]{C} that $w=w(|x|)>0$ satisfies
\begin{equation}\label{1:id}
\inte |\nabla w |^2dx  =\inte w ^2dx=\frac{1}{2}\inte w ^4dx,
\end{equation}
and note from \cite[Proposition 4.1]{GNN} that
\begin{equation}
w(x) \, , \ |\nabla w(x)| = O(|x|^{-\frac{1}{2}}e^{-|x|}) \quad
\text{as \ $|x|\to \infty$.}  \label{1:exp}
\end{equation}
Given any vector function $\mathcal{A }\in L^2_{loc}(\R^2,\R^2)$, recall also from \cite{Lieb} the following diamagnetic inequality:
\begin{equation}
|(\nabla -i\mathcal{A } )u|^2 \ge \big| \nabla |u|\big|^2 \,\ \hbox{a.e. on }\ \R^2,\ \, u\in H^1(\R^2,\mathbb{C}).
\label{2:2:1A}
\end{equation}

In this paper, we often use the following linearized operator
\begin{equation}\label{1.74}
\mathcal{L}:=-\Delta+1-w^2\quad\hbox{in}\,\ \R^2.
\end{equation}
It then obtains from   \cite{Frank} that
\begin{equation}\label{1.75}
ker\mathcal{L}=\{w\}\quad\hbox{and}\quad \langle\mathcal{L}v,v\rangle\geq 0 \ \ \mbox{for all} \ \ v\in L^2(\R^2),
\end{equation}
see also  \cite[Theorem 11.8]{Lieb} and \cite[Corollary 11.9]{Lieb}.
Further, by a standard argument (e.g. (3.45) in \cite{GLY}), there exists $\rho>0$
such that
\begin{equation}\label{1.76}
\langle\mathcal{L}v,v\rangle\geq \rho \|v\|^2_{H^1(\R^2)}\ \ \mbox{for all} \ \ v\in \mathcal{S},
\end{equation}
where the space $ \mathcal{S}$ is defined as
\[
\mathcal{S}:=\Big\{v\in H^1(\R^2,\R):\ \inte w(x)v(x)dx=0\Big\}.
\]
 Denote the linearized operator $\mathcal{N}$ by
\begin{equation}\label{1.89}
\mathcal{N}:=-\Delta+1-3w^2\quad\hbox{in}\,\ \R^2.
\end{equation}
It then follows from \cite{Frank,K,NT} that
\begin{equation}\label{1.90}
\hbox{ker }\mathcal{N}=\hbox{span}\Big\{\frac{\partial w}{\partial x_1},\ \frac{\partial w}{\partial x_2}\Big\}.
\end{equation}

Here and in the sequel, we always denote $u_a$ to be a complex-valued minimizer of $e_F(a)$, where the rotating speed $\Omega\in (0, \Omega^*)$ is fixed.
By the variational theory, there exists a Lagrange multiplier $\mu_a\in \R$ satisfying
\begin{equation}\label{1.2}
\mu_a=e_F(a)-\frac{a}{2}\inte |u_a|^4dx\end{equation}
such that $u_a$ solves the following Euler-Lagrange equation:
\begin{equation}\label{1.1}
-\Delta u_a+V(x)u_a+i\Omega (x^{\perp}\cdot \nabla u_a)=\mu_au_a+a|u_a|^2u_a\,\ \,\hbox{in}\,\ \R^2.				
\end{equation}
Under the assumptions of Theorem \ref{thm1.2}, we also define
\begin{equation}\label{2.1}
\alp_a:=\frac{(a^*-a)^{\frac{1}{2+p}}}{\lambda}>0,\ \,\mbox{where}   \,\ 1<p\le 2,
\end{equation}
and
\begin{equation}\label{def:li}
\lambda   =\arraycolsep=1.5pt\left\{\begin{array}{lll}
\Big[\displaystyle\frac{p}{2}  \inte h(x+y_0)w^2(x)dx\Big]
^{\frac{1}{2+p}}, \quad   &\mbox{if} & \quad 1<p<2;\\[4mm]
\displaystyle\Big[   \inte \Big(h(x+y_0)+\frac{\Omega ^2 }{4}|x|^2\Big)w^2(x)dx  \Big]
^{\frac{1}{4}}, \quad   &\mbox{if}& \quad p=2,\\[4mm]
\end{array}\right.
\end{equation}
where $y_0\in\R^2$ denotes the unique non-degenerate critical point of $H(y)$, see Theorem \ref{thm1.2}.
Setting
\begin{equation}\label{eea}
\varepsilon_{a}:=\Big(\inte|\nabla u_{a}|^{2}dx\Big)^{-\frac{1}{2}}>0,
\end{equation}
we now define
\begin{equation}\label{3.9}
w_a(x):=\varepsilon_a u_a\big(\varepsilon_a x+x_a\big)e^{-i (\frac{\varepsilon _a \Omega}{2} x\cdot x_{a}^{\bot}-\theta_a)},
\end{equation}
where $x_a$ is a global maximal point of $|u_a(x)|$ and $\theta_a\in [0,2\pi)$ is a proper constant.
Using above notations, since the proof of the following lemma is similar to those of \cite[Section 3]{GLY}, we omit the details of the proof for simplicity.

\begin{lem}\label{lem3.3}
Under the assumptions of Theorem \ref{thm1.2}, let $u_{a}$ be a minimizer of $e_F(a)$.
Then we have
\begin{enumerate}
\item [(i).] The parameter $\varepsilon_{a}>0$ satisfies
\begin{equation}\label{3:DDE}
\varepsilon_a=\alp_a+o(\alp_a)>0\ \ \mbox{and}\ \ \mu_a\varepsilon ^2_{a}\to -1\,\ \hbox{as }\ a\nearrow a^*,
\end{equation}
where $\mu_a\in \R$ is the Lagrange multiplier of (\ref{1.1}).
\item [(ii).]  The function $w_{a}(x)$, defined in (\ref{3.9})  for some suitable constant $\theta_a\in [0,2\pi)$, satisfies
\begin{equation}\label{1.13}
\lim_{a\nearrow a^*}w_{a}(x)=\frac{w(x) }{\sqrt{a^*}}\ \ \hbox{strongly in }\ L^\infty (\R^2, \mathbb{C})\cap H^1 (\R^2, \mathbb{C})\,\ \hbox{as }\ a\nearrow a^*,
\end{equation}
and there exists a constant $C>0$, independent of $0<a<a^*$, such that $w_{a}(x)$ satisfies
\begin{equation}\label{new4}
|w_{a}(x)|\leq Ce^{-\frac{2}{3}|x|}\quad \hbox{uniformly in}\ \,  \R^2 \,\  \hbox{as }\ a\nearrow a^*.
\end{equation}
\item [(iii).]   The global maximal point  $x_a$ of $|u_a|$ must be unique as  $a\nearrow a^*$, and $x_a$  satisfies
\begin{equation}\label{5:DDE}
\lim_{a\nearrow a^*}\frac{x_a}{\varepsilon_a}=\lim_{a\nearrow a^*}\frac{x_a}{\alp_a}=y_0,
\end{equation}
where $y_0\in\R^2$ denotes the unique non-degenerate critical point of $H(y)$.
\end{enumerate}
\end{lem}

\vskip 0.05truein

Applying Lemma  \ref{lem3.3}, we now consider
\begin{equation}\label{1.11}
v_{a}(x):=\alp_a u_{a}( \alp_a(x +y_0) )e^{-i\, \big(\frac{\Omega\alp_a^2}{2}x\cdot y_0^\perp\big)}e^{i\varphi_a}=R_a(x)+iI_a(x),
\end{equation}
where $R_a(x)$ and $I_a(x)$ denote the real and imaginary parts of $v_{a}(x)$, respectively, and the constant phase $\varphi_a\in [0,2\pi)$ is chosen such that
\begin{equation}\label{1.9}
\Big\| v_{a}-\frac{w}{\sqrt{a^*}}\Big\|_{L^2(\R^2)}=\min_{\theta\in [0,2\pi)}\Big\|e^{i\theta} v_{a}-\frac{w}{\sqrt{a^*}}\Big\|_{L^2(\R^2)}.
\end{equation}
This gives the following orthogonality condition on $I_a(x)$, which  plays an essential role in proving Theorem \ref{thm1.2}:
\begin{equation}\label{1.12}
\inte w(x)I_a(x)dx=0.
\end{equation}
Based on Lemma  \ref{lem3.3}, we next derive the following $L^\infty-$uniform estimates of $v_a$:

\begin{prop}\label{5:prop}  Under the assumptions of Theorem \ref{thm1.2},
assume $v_a(x)$ and $\varphi_a\in [0,2\pi)$ are defined by \eqref{1.11} and \eqref{1.9}, respectively. Then we have
\begin{enumerate}
\item [(i).] The function $v_{a}(x)$ satisfies
\begin{equation}
\lim_{a\nearrow a^*}v_{a}(x)=\frac{w(x) }{\sqrt{a^*}}\ \ \hbox{strongly in $L^\infty (\R^2, \mathbb{C})\cap H^1 (\R^2, \mathbb{C})$. }
\label{I:con:a}
\end{equation}
\item [(ii).] There exists a constant $C>0$ such that
\begin{equation}\label{1.19}
|v_{a}(x)|\leq Ce^{-\frac 23|x|},\ \ |\nabla v_{a}(x)|\leq Ce^{-\frac 12|x|}\ \,\hbox{in}\ \,\R^2\ \,  \hbox{as }\, \, a\nearrow a^*.	
\end{equation}
\end{enumerate}
\end{prop}

\noindent{\bf Proof.} 1.
We first derive from \eqref{3:DDE}, \eqref{1.13} and \eqref{5:DDE} that
\begin{equation}\label{1.14}
\lim_{a\nearrow a^*}\alp_a u_{a}( \alp_a(x +y_0) )e^{-i\, \big(\frac{\Omega\alp_a^2}{2}x\cdot y_0^\perp\big)}e^{i\theta_a}=\frac{w(x) }{\sqrt{a^*}}\,\ \,\hbox{strongly in $L^\infty (\R^2, \mathbb{C})\cap H^1 (\R^2, \mathbb{C})$,}
\end{equation}
where $\theta_a\in [0,2\pi)$ is defined by
\eqref{3.9} such that \eqref{1.13} holds.
We then obtain from \eqref{1.14} that
\[
\begin{aligned}
&\quad\lim_{a\nearrow a^*}\big\|\alp_a u_{a}( \alp_a(x +y_0) )e^{-i\, \big(\frac{\Omega\alp_a^2}{2}x\cdot y_0^\perp\big)}\big(e^{i\theta_a}-e^{i\varphi_a}\big)\big\|_{L^2(\R^2)}\\
&\leq \lim_{a\nearrow a^*}\Big\|\alp_a u_{a}( \alp_a(x +y_0) )e^{-i\, \big(\frac{\Omega\alp_a^2}{2}x\cdot y_0^\perp\big)}e^{i\theta_a}-\frac{w}{\sqrt{a^*}}\Big\|_{L^2(\R^2)}\\
&~~ +\lim_{a\nearrow a^*}\Big\|\alp_a u_{a}( \alp_a(x +y_0) )e^{-i\, \big(\frac{\Omega\alp_a^2}{2}x\cdot y_0^\perp\big)}e^{i\varphi_a}-\frac{w}{\sqrt{a^*}}\Big\|_{L^2(\R^2)}\\
&\leq \lim_{a\nearrow a^*}2\Big\|\alp_a u_{a}( \alp_a(x +y_0) )e^{-i\, \big(\frac{\Omega\alp_a^2}{2}x\cdot y_0^\perp\big)}e^{i\theta_a}-\frac{w}{\sqrt{a^*}}\Big\|_{L^2(\R^2)}=0,
\end{aligned}
\]
where we have used \eqref{1.9} in the last inequality.
The above estimate implies that
\begin{equation}\label{1.33}
\lim_{a\nearrow a^*}(\varphi_a-\theta_a)=0,
\end{equation}
and Proposition \ref{5:prop} (i) is therefore proved in view of \eqref{1.14}.

2. Note from \eqref{1.1} and \eqref{1.11} that $v_a(x)$ satisfies
\begin{equation}\label{1.40}
\begin{split}
&-\Delta v_a+i\alp_a^2\Omega \big(x^{\perp}\cdot\nabla v_a\big)+\Big[\frac{\alp_a^4\Omega^2|x|^2}{4}+\alp_a^2V_{\Omega}(\alp_a(x
+y_0))\Big]v_a\\
=&\alp_a^2\mu_{a}v_a+a|v_a|^2v_a\quad\hbox{in }\,\,\R^2.
\end{split}
\end{equation}
Similar to (\ref{2.21M}) below, by the comparison principle one can derive from   \eqref{I:con:a} and (\ref{1.40}) that as $a\nearrow a^*$,
\begin{equation}\label{1.44}
|v_{a}(x)|\leq Ce^{-\frac 23|x|}\quad\hbox{in}\,\ \R^2.
\end{equation}
It remains to prove the uniformly exponential decay of $|\nabla v_a|$ as $a\nearrow a^*$. Towards this aim, denoting $\partial_j v_a(x):=\frac{\partial v_a(x)}{\partial x_j}$ ($j=1,\,2$), we follow from (\ref{1.40}) that for $j,\, l=1,\,2$,
\[
\begin{aligned}
&-\Delta \partial_j v_a+(-1)^{j+1}\f^2\o i\partial_l v_a+i\f^2\o(x^{\perp}\cdot\nabla \partial_j v_a)\\
&+\Big(\frac{\f^4\o^2}{4}|x|^2+\f^2V_{\Omega}(\f (x+y_0))-\f^2\mu _a-a|v_a|^2\Big)\partial_j v_a\\
&+\Big[\frac{\f^4\o^2}{2}x_j+\f^2\frac{\partial V_{\o}(\f (x+y_0))}{\partial x_j}-2a(\partial_j v_a,v_a)\Big]v_a=0\ \ \hbox{in}\ \, \R^2,\ \ l\not =j,
\end{aligned}
\]
where $(f,g)=Re(f\cdot \bar g)$ denotes the real part of $f\cdot \bar g$.
We then get that for $l\not =j$,
\begin{equation}
\begin{aligned}
&-\frac{1}{2}\Delta |\partial_j v_a|^2+|\nabla \partial_j v_a|^2+(-1)^{j+1}\f^2\o(i\partial_l v_a,\partial_j v_a)-\f^2\o x^{\perp}\cdot(i\partial_j v_a,\nabla \partial_j v_a)\\
&+\Big[\frac{\f^4\o^2}{4}|x|^2+\f^2V_{\o}(\f(x+y_0))-\f^2\mu _a-a|v_a|^2\Big]|\partial_j v_a|^2\\
&+\Big(\frac{\f^4\o^2}{2}x_j+\f^2\frac{\partial V_{\o}(\f (x+y_0))}{\partial x_j}\Big)(v_a,\partial_j v_a)-2a(\partial_j v_a,v_a)^2=0\ \ \hbox{in}\ \, \R^2,\ \ j,\,l=1,\, 2.
\end{aligned}
\end{equation}
By the inequality (\ref{2:2:1A}), we have
\begin{equation}\label{4A:B1}
|\nabla \partial_j v_a|^2-\f^2\o x^{\perp}\cdot (i\partial_j v_a,\nabla \partial_j v_a)+\frac{\f^4\o ^2}{4}|x|^2|\partial_j v_a|^2\geq 0\quad\hbox{in\ $\R^2$}, \ \ j=1,\, 2.
\end{equation}
Under the assumption $(V)$, we also obtain from (\ref{1.44}) that for $j=1,2,$
\[
\begin{aligned}
&\Big|\Big(\frac{\f^4\o^2}{2}x_j+\f^2\frac{\partial V_{\o}(\f(x+y_0))}{\partial x_j}\Big)(v_a,\partial_j v_a)\Big|\\
\leq & \frac{\f^4\Omega^2}{2}\Big(\frac{|x_j|^2|v_a|^2}{2}+\frac{|\partial_j v_a|^2}{2}\Big)+\f^2\Big(\frac{|\frac{\partial V_{\o}(\f(x+y_0))}{\partial x_j}|^2|v_{a}|^2}{2}+\frac{|\partial_j v_a|^2}{2}\Big)\\
\leq & C\f^{2}e^{-|x|} +C\f^2|\partial_j v_a|^2 \quad \mbox{in}\ \ \R^2,
\end{aligned}\]
and
\[
\begin{aligned}
&(-\f^2\mu _a-a|v_a|^2)\big(|\partial_1v_a|^2+|\partial_2v_a|^2\big)-2a\big[(\partial_1v_a,v_a)^2+(\partial_2v_a,v_a)^2\big]\\
&+\f^2\o(i\partial_2v_a,\partial_1v_a)-\f^2\o(i\partial_1v_a,\partial_2v_a)\\
\geq & \Big(\frac{3}{4}-3a|v_a|^2\Big)\big(|\partial_1v_a|^2+|\partial_2v_a|^2\big)\quad \mbox{in}\ \ \R^2,
\end{aligned}
\]
since $-\alpha_a^2\mu_a\to 1$ and $\alpha_a^2 \to 0$ as $a\nearrow a^*$.
Combining all above estimates, we get that as $a\nearrow a^*$,
\begin{equation}\label{4:3C}
\Big(-\frac{1}{2}\Delta+\frac{2}{3}-3a|v_a|^2\Big) |\nabla v_a|^2\leq Ce^{-|x|}  \quad \mbox{in}\ \ \R^2.
\end{equation}
Because $v_a$ is bounded uniformly in $H^1(\R^2)$, one can deduce from De Giorgi-Nash-Moser theory \cite[Theorem 4.1]{HL} that
\begin{equation}\label{1.20}
\max_{x\in B_1(\xi)}|\nabla v_a(x)|^2\leq C\Big(\int_{B_2(\xi)}|\nabla v_a(x)|^2dx+\|e^{-|x|}\|_{L^2(B_2(\xi))}\Big).
\end{equation}
Therefore, we derive from \eqref{I:con:a} and \eqref{1.20} that
\begin{equation}\label{1.21}
|\nabla v_a(x)|^2\leq C\quad\hbox{and}\quad\lim_{|x|\to\infty }|\nabla v_a(x)|^2=0\ \ \hbox{as\ \ $a\nearrow a^*$.}
\end{equation}
Substituting \eqref{1.21} and \eqref{1.44} into \eqref{4:3C} yields that as $a\nearrow a^*$,
\begin{equation}\label{2.21}
\Big(-\Delta+\frac{4}{3}\Big) |\nabla v_a(x)|^2\leq Ce^{-|x|}  \quad \mbox{in}\ \ \R^2.
\end{equation}
Since $3Ce^{-|x|}>0$ is a supersolution of \eqref{2.21} in the sense that
$$\Big(-\Delta+\frac{4}{3}\Big) (3Ce^{-|x|})\geq Ce^{-|x|}
\quad\hbox{in} \,\ \R^2,$$
we have
$$
\Big(-\Delta+\frac{4}{3}\Big)\big(|\nabla v_a(x)|^2-3Ce^{-|x|}\big)\leq 0\quad\hbox{in} \,\ \R^2\ \ \hbox{as\ \ $a\nearrow a^*$.}
$$
Note from (\ref{1.21}) that $\lim_{|x|\to\infty }(|\nabla v_a(x)|^2-3Ce^{-|x|})=0$ as $a\nearrow a^*$.
We then deduce from the comparison principle that as $a\nearrow a^*$,
\begin{equation}\label{2.21M}
|\nabla v_a(x)|^2\le 3Ce^{-|x|}  \quad\hbox{in} \,\ \R^2,
\end{equation}
and the gradient estimate of (\ref{1.19}) hence holds. \qed

\vskip 0.05truein
Employing the non-degeneracy \eqref{1.76} of the linearized operator $\mathcal{L}$, we next establish the following refined  $L^\infty-$uniform estimates of the imaginary part $I_a$ for $v_a$:

\begin{lem}\label{prop-2}
Under the assumptions of Theorem \ref{thm1.2}, let $I_a$ be defined by \eqref{1.11}. Then we have
\begin{equation}\label{1.8}
|I_a(x)|\le C_{21}(\alp_a)e^{-\frac{1}{4}|x|},\ \ |\nabla I_a(x)|\le C_{22}(\alp_a)e^{-\frac{1}{8}|x|}\quad\hbox{uniformly in}\,\ \R^2\ \ \hbox{as}\,\ a\nearrow a^*,
\end{equation}  	
where the constants $C_{21}(\alp_a)>0$ and $C_{22}(\alp_a)>0$ satisfy
\begin{equation}\label{1.8M}
C_{21}(\alp_a)=o(\alp^2_a)\,\ \mbox{and} \ \ C_{22}(\alp_a)=o(\alp^2_a) \ \ \hbox{as}\,\ a\nearrow a^*.
\end{equation}
\end{lem}

\noindent{\bf Proof.} Denote the operator
\begin{equation}\label{1.48}
\mathcal{L}_a:=-\Delta +\frac{\alp_a^4\Omega^2|x|^2}{4}+\alp_a^2V_{\Omega}(\f(x+y_0))-\alp_a^2\mu_{a}-a|v_a|^2.
\end{equation}
Following \eqref{1.12} and \eqref{1.40}, the imaginary part $I_a$ of $v_a$ satisfies
\begin{equation}\label{1.15}
\mathcal{L}_aI_a(x)=-\alp_a^2\Omega\big(x^{\perp}\cdot \nabla R_a\big)\quad\hbox{in}\,\ \R^2,\quad \inte I_a(x)w(x)dx= 0.
\end{equation}
We first claim that as $a\nearrow a^*$,
\begin{equation}\label{1.16}
\big|\alp_a^2\Omega \big(x^{\perp}\cdot\nabla R_a\big)\big|\le C(\alp_a)e^{-\frac 14|x|}\quad\hbox{uniformly\,\ in\,\ $\R^2$},
\end{equation}
where  the constant $C(\alp_a)>0$ satisfies $C(\alp_a)=o(\alp^2_a)$ as $a\nearrow a^*$.
Actually, since $(x^\perp\cdot \nabla w)\equiv0$ in $\R^2$,
 for any fixed large $R>0$, one can obtain from \eqref{1:exp} and \eqref{1.19} that
\begin{equation}\label{1.47}
\Big| x^{\perp}\cdot\nabla \Big(R_a-\frac{w}{\sqrt{a^*}}\Big)\Big|\leq Ce^{-\frac{R}{8}}e^{-\frac{1}{4}|x|}\quad\hbox{in}\,\ \R^2/B_R(0).
\end{equation}
On the other hand, we get from \eqref{1.19} that $v_a$ and $ (x^{\perp}\cdot\nabla v_a)$ are bounded uniformly in $L^{\infty}\big(B_{R+1}(0)\big)$. Applying the $L^p$ estimate (cf.\cite[Theorem 9.11]{GT}) to \eqref{1.40} yields that $v_a$ is bounded uniformly in $W^{2,q}\big(B_R(0)\big)$ for any $q>2$. Since the embedding $W^{2,q}\big(B_R(0)\big)\hookrightarrow C^1\big(B_R(0)\big)$ is compact (cf.\cite[Theorem 7.26]{GT}), there exist a subsequence $\{v_{a_k}\}$ of $\{v_a\}$ and $w_0(x)\in H^1\big(B_R(0)\big)$ such that
\[
v_{a_k}(x)\to w_0(x)\quad\hbox{uniformly in\,\  $C^1\big(B_R(0)\big)$\,\ as\,\ $a_k\nearrow a^*$.}
\]
By the convergence \eqref{I:con:a} and the uniqueness of $w$, we conclude that $w_0(x)\equiv \frac{w(x)}{\sqrt{a^*}}$ in $B_R(0)$ and the above convergence hence holds for the whole sequence, i.e.,
\begin{equation}\label{1.46}
v_{a}(x)\to \frac{w}{\sqrt{a^*}}\quad\hbox{uniformly\,\ in\,\  $C^1\big(B_R(0)\big)$\,\ as\,\ $a\nearrow a^*$.}
\end{equation}
Since $R>0$ is arbitrary, we deduce from \eqref{1.47} and \eqref{1.46} that the claim \eqref{1.16} holds true.


We next follow \eqref{1.16} to prove \eqref{1.8} and  \eqref{1.8M}.   Multiplying \eqref{1.15} by $I_a$ and integrating over $\R^2$, we get from (\ref{1.16}) that
\begin{equation}\label{1.49}
\inte (\mathcal{L}_aI_a)I_adx=-\alp^2_a\Omega \inte(x^\bot\cdot\nabla {R}_a)I_adx=o(\alp_a^2)\|I_a\|_{L^2(\R^2)}\ \ \mbox{as}\ \  a\nearrow a^*.
\end{equation}
Following \eqref{1.76} and \eqref{I:con:a}, since $\inte I_awdx= 0,$  we also have
\begin{equation}\label{1.50}
\begin{split}
\inte (\mathcal{L}_aI_a)I_adx&\ge  \inte \Big\{(\mathcal{L}I_a)I_a- \big(1+\alp^2_a\mu _a\big)I_a^2- \big(a|v_a|^2-w^2\big)I_a^2\Big\}dx\\&
= \inte (\mathcal{L}I_a)I_adx+o(1)\inte I_a^2dx\geq \frac{\rho}{2}\|I_a\|^2_{H^1(\R^2)}\ \ \mbox{as}\ \  a\nearrow a^*,
\end{split}\end{equation}
where the constant $\rho>0$, independent of $0<a<a^*$, is given by \eqref{1.76}.
Hence, we obtain from \eqref{1.49} and \eqref{1.50} that
\begin{equation}\label{1.51}
\|I_a\|_{H^1(\R^2)}=o(\alp_a^2)\ \ \mbox{as}\ \  a\nearrow a^*.
\end{equation}
On the other hand, we derive from \eqref{1.15} that $|I_a|^2$ satisfies
\begin{equation*}
\begin{split}
&\Big[-\frac 12\Delta +\Big(\frac{\alp_a^4\Omega^2}{4}|x|^2+\alp_a^2V_\Omega (\f(x+y_0))-\mu _a\alp_a^2-a |v_{a}|^2\Big) \Big]|I_{a}|^2+|\nabla I_a|^2\\
= &-\alp^2_a\Omega (x^\bot\cdot\nabla {R}_{a})I_a\ \ \hbox{in}\,\ \R^2,
\end{split}
\end{equation*}
which implies that
\begin{equation}\label{1.52}
-\frac 12\Delta|I_{a}|^2-\mu _a\alp_a^2|I_{a}|^2-a |v_{a}|^2|I_{a}|^2\leq -\alp^2_a\Omega (x^\bot\cdot\nabla {R}_{a})I_a\ \ \hbox{in}\,\ \R^2.
\end{equation}
Following De Giorgi-Nash-Moser theory again (cf. \cite[Theorem 4.1]{HL}), it follows from (\ref{1.52}) that for any $\xi\in \R^2$,
\begin{equation}\label{1.53}
\sup_{x\in B_1(\xi)}|I_a(x)|^2\leq C\Big(\|I_a\|^2_{L^2(B_2(\xi))}+\|\alp^2_a\Omega (x^\bot\cdot\nabla {R}_{a})I_a\|_{L^2(B_2(\xi))}\Big).
\end{equation}
Applying Proposition \ref{5:prop}, we then deduce from \eqref{1.16}, \eqref{1.51} and \eqref{1.53} that
\begin{equation}\label{1.54}
\|I_a\|_{L^\infty(\R^2)}=o(\alp_a^2)\ \ \mbox{as}\ \  a\nearrow a^*,
\end{equation}
and hence
\begin{equation}\label{1.55}
\Big|a|v_a|^2I_a-\alp^2_a\Omega (x^\bot\cdot\nabla {R}_{a})\Big|\le C_0(\alp_a)e^{-\frac{1}{4}|x|}\quad\hbox{uniformly\,\ in\,\ $\R^2$}\ \ \mbox{as}\ \  a\nearrow a^*,
\end{equation}
where  the constant $C_0(\alp_a)>0$ satisfies $C_0(\alp_a)=o(\alp^2_a)$ as $a\nearrow a^*$.
Applying the comparison principle to \eqref{1.52}, we further get from \eqref{1.54} and \eqref{1.55} that
\[
|I_a(x)|\le C_1(\alp_a)e^{-\frac{1}{4}|x|}\quad\hbox{uniformly\,\ in\,\ $\R^2$}\ \ \mbox{as}\ \  a\nearrow a^*,
\]
where  the constant $C_1(\alp_a)>0$ satisfies $C_1(\alp_a)=o(\alp^2_a)$ as $a\nearrow a^*$.
Moreover, applying gradient estimates (see (3.15) in \cite{GT}) to the equation \eqref{1.15}, we conclude from above that the gradient estimate of \eqref{1.8} and  \eqref{1.8M} holds true, which therefore completes the proof of Lemma \ref{prop-2}. \qed

\section{Analysis of the linearized problem \eqref{1.27}}
Following the $L^\infty-$uniform estimates of previous section, this section is concerned with the analysis of the linearized problem (\ref{1.27}) defined below. In order to prove Theorem \ref{thm1.2}, by contradiction we first suppose that up to a constant phase,  there exist two different minimizers $u_{1,a}$ and $u_{2,a}$ of $e_F(a)$ as $a\nearrow a^*$, in the sense that $u_{1,a}\not\equiv u_{2,a}e^{i\theta}$ for any constant phase $\theta=\theta (a)\in[0,2\pi)$.

Recall that $\alp_a>0$ is defined by (\ref{2.1}), and $y_0\in\R^2$ denotes the unique non-degenerate critical point of $H(y)$. We then define for $j=1,\,2,$
\begin{equation}\label{3.21}
\tilde u_{j,a}(x):=\alp_au_{j,a}(\f(x+y_0))e^{-i(\frac{\alp_a^2\Omega }{2} x\cdot y_0^{\perp})}e^{i\varphi_{j,a}}=R_{j,a}(x)+iI_{j,a}(x),
\end{equation}
where  $R_{j,a}(x)$ and $I_{j,a}(x)$ denote the real and imaginary parts of $\tilde u_{j,a}(x)$, respectively, and  the constant phase $\varphi_{j,a}\in [0,2\pi)$ can be chosen properly such that
\begin{equation}\label{1.30}
\inte w(x)I_{j,a}(x)dx=0,\,\ j=1,2.
\end{equation}
We remark that (\ref{1.30}) makes sense in view of (\ref{1.9}) and the following fact: it yields from Proposition \ref{5:prop} that  $\tilde u_{j,a}(x)$ satisfies
for $j=1,2,$
\begin{equation}
\lim_{a\nearrow a^*}\tilde u_{j,a}(x)=\frac{w(x) }{\sqrt{a^*}}\,\ \hbox{strongly\,\ in\,\ $L^\infty (\R^2, \mathbb{C})\cap H^1 (\R^2, \mathbb{C})$.}
\label{1.69}
\end{equation}
Note from \eqref{1.1} and \eqref{3.21} that $\tilde u_{j,a}(x)$ satisfies the equation
\begin{equation}\label{1.67}
\begin{split}
&\quad-\Delta \tilde u_{j,a}(x)+i\alp_a^2\Omega (x^{\perp}\cdot \nabla \tilde u_{j,a})+\Big[\frac{\alp_a^4\Omega^2|x|^2}{4}+\alp_a^2V_{\Omega}(\f(x+y_0))\Big]\tilde u_{j,a}(x)\\
&=\alp_a^2\mu_{j,a}\tilde u_{j,a}(x)+a|\tilde u_{j,a}|^2\tilde u_{j,a}(x)\quad\hbox{in}\,\ \R^2,\,\ j=1,2.
\end{split}
\end{equation}
Applying \eqref{1.30}-\eqref{1.67}, we have the following estimates:

\begin{lem}\label{prop-1}
Under the assumptions of Theorem~\ref{thm1.2}, let $I_{j,a}(x)$ be defined by \eqref{3.21} for $j=1,2$. Then $I_{j,a}(x)$ satisfies	for $j=1,\,2$,
\begin{equation}\label{3.23}
|I_{j,a}(x)|\le C_{j1}(\alp_a)e^{-\frac 14 |x|}\ \, \mbox{and}\ \ |\nabla I_{j,a}(x)|\le C_{j2}(\alp_a)e^{-\frac 18 |x|}\,\ \hbox{uniformly\,\ in}\,\ \R^2\ \, \mbox{as}\,\  a\nearrow a^*,
\end{equation}
where the constants $C_{j1}(\alp_a)>0$ and $C_{j2}(\alp_a)>0$ satisfy
\begin{equation}\label{3.23M}
C_{j1}(\alp_a)=o(\alp^2_a)\,\ \mbox{and} \ \ C_{j2}(\alp_a)=o(\alp^2_a) \ \ \hbox{as}\,\ a\nearrow a^*.
\end{equation}
\end{lem}


\vskip 0.05truein

Since Lemma \ref{prop-1} can be established in the similar approach of proving Lemma \ref{prop-2}, we omit the detailed proof for simplicity.
We next define the following difference function:
\begin{equation}\label{1.24}
\eta_a(x):=\frac{\tilde u_{2,a}(x)-\tilde u_{1,a}(x)}{\|\tilde u_{2,a}(x)-\tilde u_{1,a}(x)\|_{L^{\infty}(\R^2)}}=\eta_{1,a}(x)+i\eta_{2,a}(x),
\end{equation}
where $\eta_{1,a}(x)$ and $\eta_{2,a}(x)$ denote the real and  imaginary parts of $\eta_a(x)$, respectively. By the definition of $\eta_a(x)$, we deduce from \eqref{1.67} that $\eta_a(x)$ satisfies
\begin{equation}\label{1.27}
\begin{split}
&~~-\Delta \eta_a+i\alp_a^2\Omega(x^{\perp}\cdot \nabla \eta_a)+\alp_a^2V_{\Omega}(\f(x+y_0))\eta_a+\frac{\alp_a^4\Omega^2|x|^2}{4}\eta_a\\
&=\alp_a^2\mu_{1,a}\eta_a-\frac{\tilde u_{2,a}}{2}\inte \tilde f_a(|\tilde u_{1,a}|^2+|\tilde u_{2,a}|^2)dx+\tilde f_a\tilde u_{1,a}+a|\tilde u_{2,a}|^2\eta_a\,\ \hbox{in}\,\ \R^2,
\end{split}
\end{equation}
where $\tilde f_a$ is defined by
\begin{equation}\label{1.62}
\begin{split}
\tilde f_a(x):&=a\frac{|\tilde u_{2,a}(x)|^2-|\tilde u_{1,a}(x)|^2}{\|\tilde u_{2,a}-\tilde u_{1,a}\|_{L^\infty(\R^2)}}=a\Big[\eta_{1,a}\big(R_{1,a}+R_{2,a}\big)+\eta_{2,a}\big(I_{1,a}+I_{2,a}\big)\Big].
\end{split}
\end{equation}
We next study the estimates of $\eta_a(x)$ as $a\nearrow a^*$.

\begin{lem}\label{1.25}
Suppose $\eta_a(x)$ is defined by \eqref{1.24}, then there exists a constant $C>0$, independent of $0<a<a^*$, such that
\begin{equation}\label{1.26}
|\eta_a(x)|\leq Ce^{-\frac 23|x|}\ \ \hbox{and}\quad |\nabla\eta_a(x)|\leq Ce^{-\frac 12|x|}\,\ \hbox{uniformly in\,\ $\R^2$}\, \ \hbox{as}\, \ a\nearrow a^*.
\end{equation}
\end{lem}

\noindent\textbf{Proof.} We first address the estimate of $\eta_a(x)$ as $a\nearrow a^*$.
Following \eqref{1.27}, we have
\begin{equation}\label{1.80}
\begin{split}
&~~-\frac 12\Delta |\eta_a|^2+|\nabla \eta_a|^2-\alp_a^2\Omega x^{\perp}\cdot(i\eta_a, \nabla \eta_a)+\Big[\alp_a^2V_{\Omega}(\f(x+y_0))+\frac{\alp_a^4\Omega^2|x|^2}{4}\Big]|\eta_a|^2\\
&=\alp_a^2\mu_{1,a}|\eta_a|^2-\frac{(\tilde u_{2,a},\eta_a)}{2}\inte \tilde f_a(|\tilde u_{1,a}|^2+|\tilde u_{2,a}|^2)dx\\
&\quad+(\tilde f_a\tilde u_{1,a},\eta_a)+a|\tilde u_{2,a}|^2|\eta_a|^2\quad\hbox{in}\,\ \R^2.
\end{split}
\end{equation}
Using the diamagnetic inequality (\ref{2:2:1A}), we obtain that
\begin{equation}\label{1.78}
|\nabla \eta_a|^2-\alp_a^2\Omega x^{\perp}\cdot(i\eta_a, \nabla \eta_a)+\frac{\alp_a^4\Omega^2|x|^2}{4}|\eta_a|^2\geq 0\quad\hbox{in}\,\ \R^2.
\end{equation}
By the definition of $\tilde{f}_a$, we also get from Proposition \ref{5:prop} that there exists a constant $C>0$ such that
\begin{equation}\label{3.17A}
\|\tilde f_a\|_{L^\infty(\R^2)}\leq C\quad\hbox{and}\quad\Big|\inte \tilde f_a(|\tilde u_{1,a}|^2+|\tilde u_{2,a}|^2)dx\Big|\leq C.
\end{equation}
Consequently, we obtain from Proposition \ref{5:prop} that
\begin{equation}\label{1.79}
\begin{split}
&\quad\Big|-\frac{(\tilde u_{2,a},\eta_a)}{2}\inte \tilde f_a(|\tilde u_{1,a}|^2+|\tilde u_{2,a}|^2)dx+(\tilde f_a\tilde u_{1,a},\eta_a)\Big|\\
&\leq 2\delta |\eta_a|^2+C(\delta)|\tilde u_{2,a}|^2+C(\delta)|\tilde f_a|^2|\tilde u_{1,a}|^2\\
&\leq 2\delta |\eta_a|^2+C(\delta)e^{-\frac 43|x|}\quad\hbox{in}\,\ \R^2,
\end{split}
\end{equation}
where $\delta>0$ is a small constant and $C(\delta)>0$ depends only on $\delta$.
Recall that $-\alp_a^2\mu_{1,a}\to 1$ as $a\nearrow a^*$ in view of Lemma \ref{lem3.3}. Setting $\delta=\frac{1}{16}$, we then deduce from \eqref{1.80}-\eqref{1.79} that as $a\nearrow a^*$,
\[
-\frac 12\Delta |\eta_a|^2+\frac{3}{4}|\eta_a|^2\leq Ce^{-\frac 43|x|}\quad \hbox{in}\,\ \R^2,
\]
where Proposition \ref{5:prop} is used again.
By the comparison principle, we conclude from above that there exists a constant $C>0$, independent of $0<a<a^*$, such that as $a\nearrow a^*$,
\begin{equation}\label{1.81}
|\eta_a(x)|\leq Ce^{-\frac 23|x|}\quad \hbox{uniformly in}\,\ \R^2,
\end{equation}
which implies that the estimate (\ref{1.26}) of $\eta_a(x)$ holds true.

We next prove the estimate of $|\nabla \eta_{a}|$ as $a\nearrow a^*$. Actually, note from \eqref{1.27} and \eqref{3.17A} that as $a\nearrow a^*$,
\[
\begin{aligned}
\inte |\nabla\eta_a|^2dx&\le\inte \Big\{|\nabla\eta_a|^2+\Big[\alp_a^2V_{\Omega}(\f(x+y_0))+\frac{\alp_a^4\Omega^2|x|^2}{4}-\alp_a^2\mu_{1,a}\Big]|\eta_a|^2\Big\}dx\\
&=\inte\Big\{\alp_a^2\Omega\, x^{\perp}\cdot(i\eta_a,\nabla \eta_a)+(\tilde f_a\tilde u_{1,a},\eta_a)+a|\tilde u_{2,a}|^2|\eta_a|^2\Big\}dx\\
&\quad-\inte\frac{(\tilde u_{2,a},\eta_a)}{2}dx\inte \tilde f_a(|\tilde u_{1,a}|^2+|\tilde u_{2,a}|^2)dx\\
&\leq \inte \Big\{\frac 12|\nabla \eta_a|^2+\frac{\alp_a^4\Omega^2|x|^2}{2}|\eta_a|^2\Big\}dx+C\\
&\leq \frac12\inte |\nabla\eta_a|^2dx+C,
\end{aligned}
\]
where the estimate \eqref{1.81} is also used. Following the above estimate, there exists a constant $C>0$ such that as $a\nearrow a^*$,
\begin{equation}\label{1.43}
\|\nabla \eta_a(x)\|_{L^2(\R^2)}\leq C.
\end{equation}
For convenience, we now denote
\[
\partial_j\eta_{a}:=\frac{\partial\eta_{a}(x)}{\partial x_j},\,\ V_a(x):=\alp_a^2V_{\Omega}(\f(x+y_0))+\frac{\alp_a^4\Omega^2|x|^2}{4},
\]
\[
F_a(x):=-\frac{\tilde u_{2,a}}{2}\inte \tilde f_a(|\tilde u_{1,a}|^2+|\tilde u_{2,a}|^2)dx+\tilde f_a\tilde u_{1,a}.
\]
It follows from \eqref{1.27} that for $j,l=1,2$,
\[
\begin{aligned}
&~~-\Delta \partial_j\eta_a+i\alp_a^2\Omega(x^{\perp}\cdot \nabla \partial_j\eta_a)+(-1)^{j+1}i\partial_l\eta_a+V_a(x)\partial_j\eta_a+\frac{\partial V_{a}}{\partial x_j}\eta_a\\
&=\alp_a^2\mu_{1,a}\partial_j\eta_a+\frac{\partial F_{a}}{\partial x_j}+a|\tilde u_{2,a}|^2\partial_j\eta_a+a\frac{\partial |\tilde u_{2,a}|^2}{\partial x_j}\eta_a\ \ \hbox{in}\,\ \R^2,\ \, j\neq l.
\end{aligned}
\]
This further implies that for $j,l=1,2$,
\begin{equation}\label{1.84}
\begin{split}
&~~-\frac 12\Delta |\partial_j\eta_a|^2+|\nabla \partial_j\eta_a |^2-\alp_a^2\Omega x^{\perp}\cdot(i\partial_j\eta_a, \nabla \partial_j\eta_a)\\
&~~+(-1)^{j+1}(i\partial_l\eta_a,\partial_j\eta_a)+V_a(x)|\partial_j\eta_a|^2+\frac{\partial V_{a}}{\partial x_j}(\eta_a,\partial_j\eta_a)\\
&=\alp_a^2\mu_{1,a}|\partial_j\eta_a|^2+\Big(\frac{\partial F_{a}}{\partial x_j},\partial_j\eta_a\Big)+a|\tilde u_{2,a}|^2|\partial_j\eta_a|^2+a\frac{\partial |\tilde u_{2,a}|^2}{\partial x_j}(\eta_a,\partial_j\eta_a)\,\ \hbox{in}\,\ \R^2,\ \, j\neq l.
\end{split}
\end{equation}
Applying the comparison principle, we follow \eqref{1.81}--(\ref{1.84}) to finish the proof as follows:

By Cauchy's inequality, we get that
\begin{equation}\label{1.85}
\alp_a^2\Omega x^\perp\cdot(i\partial_j\eta_a,\nabla \partial_j\eta_a)\leq  \frac{\alp_a^4\Omega^2|x|^2}{4}|\partial_j\eta_a|^2+|\nabla \partial_j\eta_a|^2\quad \hbox{in}\,\ \R^2.
\end{equation}
Under the assumption $(V)$, we obtain from \eqref{1.19} and \eqref{1.81} that
\begin{equation}\label{1.86}
\begin{split}
\frac{\partial V_a}{\partial x_j}(\eta_a,\partial_j\eta_a)
&\leq \delta |\partial_j\eta_a|^2+C(\delta)\Big|\frac{\partial V_a}{\partial x_j}\Big|^2|\eta_a|^2\\
&\leq \delta |\partial_j\eta_a|^2+C(\delta)e^{-|x|}\quad \hbox{in}\,\ \R^2,
\end{split}
\end{equation}
and
\begin{equation}\label{1.87}
\begin{split}
a\frac{\partial |u_{2,a}|^2}{\partial x_j}(\eta_a,\partial_j\eta_a)&\leq  \delta |\partial_j\eta_a|^2+C(\delta)a^2\Big|\frac{\partial |\tilde u_{2,a}|^2}{\partial x_j}\Big|^2|\eta_a|^2\\
&\leq \delta |\partial_j\eta_a|^2+C(\delta)e^{-|x|}\quad \hbox{in}\,\ \R^2,
\end{split}
\end{equation}
where $\delta>0$ is a small constant and $C(\delta)>0$ depends only on $\delta$.
Similarly, one can derive from \eqref{1.62} and \eqref{3.17A} that
\begin{equation}\label{1.83}
\begin{split}
\quad\Big(\frac{\partial F_a}{\partial x_j},\partial_j\eta_a\Big)
&\leq 2\delta |\partial_j\eta_a|^2+C(\delta)\Big(\Big|\frac{\partial \tilde{u}_{2,a}}{\partial x_j}\Big|^2+ \Big|\frac{\partial \tilde{u}_{1,a}}{\partial x_j}\Big|^2\Big)+\Big(\frac{\partial \tilde f_a}{\partial x_j}\tilde{u}_{1,a},\partial_j \eta_a\Big)\\
&\leq 3\delta |\partial_j\eta_a|^2+C(\delta)e^{-|x|}+Ce^{- |x|}|\partial_j \eta_a|^2\quad\hbox{in}\,\ \R^2,
\end{split}
\end{equation}
where the estimates (\ref{1.19}) and (\ref{1.43}) are also used.
Moreover, since $-\alp_a^2\mu_{1,a}\to 1$ as $a\nearrow a^*$, we get that as $a\nearrow a^*$,
\begin{equation}\label{1.88}
\begin{split}
&~~-\alp_a^2\mu_{1,a}(|\partial_1\eta_a|^2+|\partial_2\eta_a|^2)+\alp_a^2(i\partial_2\eta_a,\partial_1\eta_a)-\alp_a^2(i\partial_1\eta_a,\partial_2\eta_a)\\
&\geq (-\alp_a^2\mu_{1,a}-\alp_a^2)(|\partial_1\eta_a|^2+|\partial_2\eta_a|^2)\geq \frac 34(|\partial_1\eta_a|^2+|\partial_2\eta_a|^2).
\end{split}
\end{equation}
Following \eqref{1.85}--\eqref{1.88} with $\delta=\frac{1}{60}$, we conclude from \eqref{1.84} that as $a\nearrow a^*$,
\begin{equation}\label{1.28}
-\frac 12\Delta |\nabla \eta_a|^2+\frac 23|\nabla \eta_a|^2-Ce^{-|x|}|\nabla \eta_a|^2\leq Ce^{-|x|}\quad\hbox{in}\,\ \R^2.
\end{equation}
Applying De Giorgi-Nash-Moser theory again, we then deduce from  \eqref{1.28} that for any $\xi\in\R^2$,
\[
\max_{x\in B_1(\xi)}|\nabla \eta_a(x)|^2\leq C\Big(\int_{B_2(\xi)}|\nabla \eta_a(x)|^2dx+\|e^{-|x|}\|_{L^2(B_2(\xi))}\Big)\  \ \hbox{as}\, \ a\nearrow a^*.
\]
Applying \eqref{1.43}, we thus obtain from above that as $a\nearrow a^*$,
\begin{equation}\label{1.82}
\|\nabla \eta_a(x)\|_{L^\infty(\R^2)}\leq C
\ \ \hbox{and}\,\ \lim_{|x|\to\infty }|\nabla \eta_a(x)|=0.
\end{equation}
Following \eqref{1.82}, we derive from \eqref{1.28} that as $a\nearrow a^*$,
\[
-\frac 12\Delta |\nabla \eta_a|^2+\frac 23|\nabla \eta_a|^2\leq Ce^{-|x|}|\nabla \eta_a|^2+Ce^{-|x|}\leq Ce^{-|x|}\quad\hbox{in}\,\ \R^2,
\]
from which the desired estimate (\ref{1.26}) of $|\nabla \eta_a|$ can be proved by the comparison principle. This completes the proof of Lemma \ref{1.25}. \qed
\vskip 0.05truein

\begin{lem}\label{1.23} Suppose $\eta_a(x)$ is defined by \eqref{1.24}. Then there exist a subsequence of $\{\eta_{a}\}$ (still denoted by $\{\eta_{a}$\}) and some constants $b_0$, $b_1$ and $b_2$ such that
\begin{equation}
\eta_a(x)=b_0\big(w+x\cdot \nabla w\big)+\sum^2_{j=1}b_j\frac{\partial w}{\partial x_j}+M_a(x)\,\ \hbox{uniformly\,\ in}\,\ \R^2\ \, \mbox{as}\,\  a\nearrow a^*,
\label{3d:5a}
\end{equation}	
where the lower order term $M_a(x)$ satisfies
\begin{equation}\label{1.41}
|M_a(x)|\le C_M (\alp_a)e^{-\frac 12|x|},\ \ |\nabla M_a(x)|\le C_M (\alp_a)e^{-\frac 14|x|}\,\ \hbox{uniformly\,\ in}\,\ \R^2\ \, \mbox{as}\,\  a\nearrow a^*
\end{equation}
for some constant $C_M (\alp_a)>0$ satisfying $C_M (\alp_a)=o(1)$ as $a\nearrow a^*.$
\end{lem}

\noindent\textbf{Proof.} We first claim that there exist a subsequence of $\{\eta_{a}\}$ (still denoted by $\{\eta_{a}\}$) and $b_0$, $b_1$ and $b_2$ such that
\begin{equation}\label{1.42}
\eta_a(x)\to b_0\big(w+x\cdot \nabla w\big)+\sum^2_{j=1}b_j\frac{\partial w}{\partial x_j}\,\ \hbox{\,\ uniformly\,\ in\,\ $C^1_{loc}(\R^2, \mathbb{C})$}\ \, \mbox{as}\,\  a\nearrow a^*.
\end{equation}
Indeed, we note from Lemma \ref{1.25} that the term $( x^\perp\cdot\nabla\eta_a)$ is bounded uniformly and decays exponentially for sufficiently large $|x|$ as $a\nearrow a^*$. By the standard elliptic regularity (cf. \cite{GT}), it then follows from   \eqref{1.27} and \eqref{3.17A} that $\eta _a\in C^{1,\alpha }_{loc}(\R^2, \mathbb{C})$ and $\|\eta _a\|_{C^{1,\alpha }_{loc}(\R^2, \mathbb{C})}\le C$ uniformly as $a\nearrow a^*$ for some $\alp \in (0,1)$. By the equation \eqref{1.27}, this implies that there exists a subsequence of $\{\eta_{a}\}$ (still denoted by $\{\eta_{a}\}$) such that
$\eta_a=\eta_{1,a}+i\eta_{2,a}\to \eta_0=\eta_1+i\eta_2\in H^1(\R^2,\mathbb{C})$ uniformly in $C^1_{loc}(\R^2, \mathbb{C})$ as $a\nearrow a^*$, where $\eta_0$ is a weak solution of
\[
-\Delta \eta_0+\eta_0-2w^2\eta_1-w^2\eta_0=-\Big(\frac{2}{a^*}\int_{\R^2}w^3\eta_1dx\Big)w\,\ \mbox{in}\,\ \R^2,
\]	
which further implies that $(\eta_1,\eta_2)$ satisfies the following system
\begin{equation}\label{1.29}
\left\{
\begin{aligned}
-\Delta \eta_1+\eta_1-3w^2\eta_1&=-\Big(\frac{2}{a^*}\int_{\R^2}w^3\eta_1dx\Big)w
\,\ \mbox{in}\,\ \R^2,\\[2mm]
-\Delta \eta_2+\eta_2-w^2\eta_2&=0
\,\ \mbox{in}\,\ \R^2.\\
\end{aligned}
\right.
\end{equation}
Notice from \eqref{1.30} and \eqref{1.24} that  $\inte w(x)\eta_{2,a}(x)dx\equiv0$, which further gives that $\inte w(x)\eta_2(x)dx=0$, in view of the fact that $\eta_{2,a}\to \eta_2$ in $C^1_{loc}(\R^2)$ as $a\nearrow a^*$. Since
$$(-\Delta+1-3w^2)(w+x\cdot\nabla w)=-2w,$$
we conclude from \eqref{1.75},  \eqref{1.90} and \eqref{1.29} that there exist some constants $b_0$, $b_1$ and $b_2$ such that
\[
\eta_1=b_0(w+x\cdot\nabla w)+\sum_{i=1}^2b_i\frac{\partial w}{\partial x_i}\quad\hbox{and}\quad \eta_2\equiv 0\quad\hbox{in}\,\ \R^2,
\]
and the claim (\ref{1.42}) is therefore proved.

On the other hand, for any fixed sufficiently large $R>0$, we derive from \eqref{1:exp} and \eqref{1.26} that as $a\nearrow a^*$,
\begin{equation}\label{1.57}
\Big|\eta_a(x)-\Big[b_0\big(w+x\cdot \nabla w\big)+\sum^2_{j=1}b_j\frac{\partial w}{\partial x_j}\Big]\Big|\leq Ce^{-\frac{1}{12} R}e^{-\frac 12 |x|}\quad \hbox{in\,\ $\R^2/B_R(0)$},
\end{equation}
and
\begin{equation}\label{1.58}
\Big|\nabla\eta_a(x)-\nabla\Big[b_0\big(w+x\cdot \nabla w\big)+\sum^2_{j=1}b_j\frac{\partial w}{\partial x_j}\Big]\Big|\leq Ce^{-\frac 14 R}e^{-\frac 14 |x|}\quad\hbox{in\,\ $\R^2/B_R(0)$}.
\end{equation}
Since $R>0$ is arbitrary, together with \eqref{1.42}, we conclude from \eqref{1.57} and \eqref{1.58} that \eqref{1.41} holds true, which completes the proof of Lemma \ref{1.23}.   \qed

Based on Lemma \ref{1.23}, we next establish the following refined estimate of $\eta_{a}$ as $a\nearrow a^*$.

\begin{prop}\label{p-3.3} Suppose $\{\eta_a\}$ is the subsequence obtained in Lemma \ref{1.23}. Then the imaginary part $\eta_{2,a}$ of $\eta_{a}$ satisfies
\begin{equation}\label{1.17}
\eta_{2,a}(x)=\frac{\alp_a^2\Omega}{2}\big(-b_1x_2+b_2x_1\big)w(x)+E_a(x)\,\ \hbox{uniformly\,\ in}\,\ \R^2\ \, \mbox{as}\,\  a\nearrow a^*,
\end{equation}
where $(x_1,x_2)=x\in\R^2$, the constants $b_1$ and $b_2$ are as in Lemma \ref{1.23}, and the lower order term $E_a(x)$ satisfies
\begin{equation}\label{1.18}
|E_a(x)|\le C_E (\alp_a)e^{-\frac 18|x|},\,\ |\nabla E_a(x)|\le C_E (\alp_a)e^{-\frac{1}{16}|x|}\,\ \hbox{uniformly\,\ in}\,\ \R^2\ \, \mbox{as}\,\  a\nearrow a^*
\end{equation}
for some constant $C_E (\alp_a)>0$ satisfying $C_E (\alp_a)=o(\alp_a^2)$ as $a\nearrow a^*.$
\end{prop}

\noindent\textbf{Proof.}
We first get from \eqref{1.27} that
\begin{equation}\label{1.61}
\begin{split}
\mathcal{L}_{2,a}\eta_{2,a}=-\alp_a^2\Omega(x^{\perp}\cdot \nabla \eta_{1,a})-\frac{I_{2,a}}{2}\inte \tilde{f}_a(|\tilde{u}_{1,a}|^2+|\tilde{u}_{2,a}|^2)dx+\tilde{f}_aI_{1,a}\,\ \hbox{in}\,\ \R^2,
\end{split}
\end{equation}
where $\tilde{f}_a$ is defined by \eqref{1.62} and  $\mathcal{L}_{j,a}$ is defined for $j=1,2$,
\begin{equation}\label{1.72}
\mathcal{L}_{j,a}:=-\Delta +\alp_a^2V_{\Omega}(\alp_a (x+y_0))+\frac{\alp_a^4\Omega^2|x|^2}{4}-\alp_a^2\mu_{j,a}-a|\tilde u_{j,a}|^2.
\end{equation}
Set
\[
G_a(x):=-\frac{I_{2,a}}{2}\inte \tilde{f}_a(|\tilde{u}_{1,a}|^2+|\tilde{u}_{2,a}|^2)dx+\tilde{f}_aI_{1,a}.
\]
We then get from \eqref{3d:5a} and (\ref{1.61}) that $\eta_{2,a}$ satisfies
\begin{equation}\label{D1*}
\mathcal{L}_{2,a}\eta_{2,a}=-\alp_a^2\Omega\Big[-b_1\frac{\partial w}{\partial x_2}+b_2\frac{\partial w}{\partial x_1}+Re(x^{\perp}\cdot \nabla  M_{a})\Big]+G_a(x)\,\ \hbox{in}\,\ \R^2,\,\ \inte \eta_{2,a}wdx= 0,
\end{equation}
where (\ref{1.30}) is used and $Re(\cdot)$ denotes the real part. Here the constants $b_1$ and $b_2$ are as in Lemma \ref{1.23}.
Note that $-\frac{1}{2}x_jw(x)$ is the unique solution of the following equation:
\begin{equation}\label{1.32}
-\Delta u+u-w^2u=\frac{\partial w}{\partial x_j},\quad\inte uwdx=0,\,\ j=1,2.
\end{equation}
Denote
\[
E_a(x):=\eta_{2,a}(x)-\frac{\alp_a^2\Omega}{2}\big(-b_1x_2+b_2x_1\big)w(x).
\]
Applying \eqref{1.32}, we thus obtain from  (\ref{D1*}) that $E_a(x)$ satisfies
\begin{equation}\label{1.59}
\mathcal{L}_{2,a}E_a(x)=\frac{\alp_a^2\Omega}{2}(\mathcal{L}-\mathcal{L}_{2,a})\big[(-b_1x_2+b_2x_1)w\big]- \alp_a^2\Omega\, Re(x^{\perp}\cdot \nabla M_{a})+G_a(x)\quad\hbox{in}\,\ \R^2,
\end{equation}
and
\begin{equation}\label{1.65}
\inte E_{a}(x)w(x)dx=0,
\end{equation}
where the operator $\mathcal{L}$ is defined by \eqref{1.74}.


We now estimate the right hand side of \eqref{1.59}.
By the definition of $\mathcal{L}_{2,a}$ in (\ref{1.72}), we have
\[
\frac{\alp_a^2\Omega}{2}\big|(\mathcal{L}_{2,a}-\mathcal{L})\big[(-b_1x_2+b_2x_1)w\big]\big|\le C(\alp_a)e^{-\frac 14|x|}\ \ \hbox{uniformly\,\ in\,\ $\R^2$}\ \, \mbox{as}\,\  a\nearrow a^*.
\]
where $C(\alp_a)>0$ satisfies
\begin{equation}\label{1.64M}
C(\alp_a)=o(\alp_a^2)\ \  \mbox{as}\,\  a\nearrow a^*.
\end{equation}
Applying Lemma \ref{1.23} for the estimate of $\nabla M_{a}$, we also obtain that
\[
\big|Re(\alp_a^2\Omega x^{\perp}\cdot \nabla M_{a})\big|\le C(\alp_a) e^{-\frac 18|x|}\,\ \hbox{uniformly\,\ in\,\ $\R^2$}\ \, \mbox{as}\,\  a\nearrow a^*,
\]
where $C(\alp_a)>0$ also satisfies (\ref{1.64M}).
Moreover, it follows from Lemma \ref{prop-1} and (\ref{3.17A}) that
\[
|G_a(x)|\le C(\alp_a)e^{-\frac 14|x|}\,\ \hbox{uniformly\,\ in\,\  $\R^2$}\ \, \mbox{as}\,\  a\nearrow a^*,
\]
where $C(\alp_a)>0$ satisfies (\ref{1.64M}) again. Following above estimates, we deduce from (\ref{1.59}) that
\begin{equation}\label{1.64}
\begin{split}
\big|\mathcal{L}_{2,a}E_a(x)\big|&=\Big|\frac{\alp_a^2\Omega}{2}(\mathcal{L}_{2,a}-\mathcal{L})\big[(-b_1x_2+b_2x_1)w\big]-Re(\alp_a^2\Omega x^{\perp}\cdot \nabla M_{a})+G_a(x)\Big|\\
&\le C(\alp_a)e^{-\frac 18|x|}\,\ \hbox{uniformly\,\ in\,\  $\R^2$}\ \, \mbox{as}\,\  a\nearrow a^*,
\end{split}
\end{equation}
where $C(\alp_a)>0$ satisfies (\ref{1.64M}). By the same argument of Lemma \ref{prop-2}, Proposition \ref{p-3.3} is therefore complete in view of (\ref{1.59}), \eqref{1.65} and \eqref{1.64}. \qed

\section{Proof of Theorem \ref{thm1.2}}
As outlined in the introduction, in this section we are ready to complete the proof of Theorem \ref{thm1.2} on the local uniqueness of minimizers for $e_F(a)$ as $a\nearrow a^*$.
\vskip 0.05truein


\noindent\textbf{Proof of Theorem \ref{thm1.2}.}  By contradiction, suppose that, up to a constant phase, there exist
two different minimizers $u_{1,a}$ and $ u_{2,a}$ of $e_F(a)$ as $a\nearrow a^*$. This means that $u_{1,a}\not\equiv u_{2,a}e^{i\theta}$ for any constant phase $\theta=\theta (a)\in[0,2\pi)$.
Recall that $\tilde u_{j,a}:=R_{j,a}(x)+iI_{j,a}(x)$ defined by \eqref{3.21} satisfies the following equation
\begin{equation}
\begin{split}
&-\Delta \tilde u_{j,a}(x) +i\alp^2_a\Omega(x^{\perp}\cdot \nabla \tilde u_{j,a})+\Big[\frac{\alp_a^4\Omega^2|x|^2}{4}+\alp_a^2V_{\Omega}\big(\alp_a(x+y_0)\big)\Big]\tilde u_{j,a}(x)\\
&\quad\quad\quad=\alp^2_a\mu_{j,a}\tilde  u_{j,a}(x) +a|\tilde u_{j,a}|^2\tilde u_{j,a}(x)\quad \text{in\,\ $\R^2$}, \,\ j=1,2,
\label{5.2:0}
\end{split}
\end{equation}
where $\mu_{j,a}\in\R$ satisfies
\begin{equation}\label{1.34}
\mu_{j,a}=e_F(a)-\frac{a}{2\alp_a^2}\inte |\tilde u_{j,a}|^4dx.
\end{equation}
Thus, the real part $R_{j,a}(x)$ of $\tilde u_{j,a}(x)$ satisfies
\begin{equation}\label{1.3}
\begin{aligned}
&-\Delta R_{j,a}(x)-\alp_a^2\Omega (x^{\perp}\cdot \nabla  I_{j,a})+\Big[\frac{\alp_a^4\Omega^2|x|^2}{4}+\alp_a^2V_{\Omega}\big(\alp_a(x+y_0)
\big)\Big] R_{j,a}(x)\\
&\quad\quad\quad=\alp_a^2\mu_{j,a} R_{j,a}(x)+a|\tilde u_{j,a}|^2 R_{j,a}(x)\quad \text{in\,\ $\R^2$}, \,\ j=1,2.
\end{aligned}
\end{equation}
Applying (\ref{1.34}), one can derive from Lemmas \ref{prop-1} and \ref{1.23} that as $a\nearrow a^*$,
\begin{equation}\label{4.7b}
\begin{split}
&\quad \frac{\alp_a^2(\mu_{2,a}-\mu_{1,a})}{\|\tilde u_{2,a}-\tilde u_{1,a}\|_{L^\infty(\R^2)}} \\
&=-\frac{a}{2}\inte \big[{\eta}_{1,a}({R}_{1,a}+{R}_{2,a})+\eta_{2,a}( I_{1,a}+ I_{2,a})\big]\cdot\big(|\tilde u_{2,a}|^2+|\tilde u_{1,a}|^2\big)dx\\
&=-\frac{2a}{(a^*)^{3/2}}\inte  \Big[b_0\big(w+x\cdot \nabla w\big)+\sum^2_{j=1}b_j\frac{\partial w}{\partial x_j}\Big]w^3dx+o(1)=b_0O(1)+o(1),
\end{split}
\end{equation}
where the constants $b_0$, $b_1$ and $b_2$ are as in \eqref{3d:5a}.

We next prove $b_i=0$ for (\ref{4.7b}), where $i=0,1,2$. Our main idea is to derive the linear equations of $b_i$ through Pohozaev identities of $R_{j,a}$ satisfying (\ref{1.3}), where $i=0,1,2$ and $j=1,2$. More precisely, we shall first construct Pohozaev identities of $R_{j,a}$, which are involved with the homogeneous potential $h(x)$ and the terms produced by the rotation. Applying Lemma \ref{prop-1} and Proposition \ref{p-3.3}, we shall prove that the terms produced by the rotation are lower order as $a\nearrow a^*$, which then yield the quantitative relationships between $\eta_1$ satisfying (\ref{1.29}) and the homogeneous potential $h(x)$. Inserting  \eqref{3d:5a} into those quantitative relationships, we further obtain the linear equations of $b_i$, from which we shall finally derive that $b_i=0$ for $i=0,1,2$  under some suitable assumptions on $h(x)$. We now carry out the above idea by the following three steps:


\vskip 0.05truein

\noindent{\em  Step 1.}
We claim that the constants $b_0$, $b_1$ and $b_2$ defined in \eqref{3d:5a} satisfy the following equations:
\begin{equation}\label{3.34}
b_0\inte \frac{\partial h(x+y_0)}{\partial x_l}(x\cdot \nabla w^2)dx+\sum_{j=1}^2b_j\inte \frac{\partial h(x+y_0)}{\partial x_l}\frac{\partial w^2}{\partial x_j}dx=0,\,\ l=1,2.
\end{equation}
The intuition behind (\ref{3.34}) is as follows: one can obtain the relationship between the homogeneous potential $h(x)$ and the coefficients $b_0$, $b_1$ and $b_2$ by constructing Pohozaev identities; Moreover, if one can further obtain that $b_0=0$ (see Step 2 below), then it implies from (\ref{3.34}) that $b_1=b_2=0$ under the non-degeneracy assumption of the critical point $y_0$ for $H(y)$.



To prove the above claim (\ref{3.34}), we first multiply the equation (\ref{1.3}) by $\frac{\partial R_{j,a}(x)}{\partial x_l}$ to get that
\begin{equation}\label{3.1}
\begin{split}
&\quad-\Delta R_{j,a}(x)\frac{\partial R_{j,a}(x)}{\partial x_l}-\alp_a^2\Omega (x^{\perp}\cdot \nabla I_{j,a}) \frac{\partial R_{j,a}(x)}{\partial x_l} \\
&\quad+\frac{1}{2}\Big[\frac{\alp_a^4\Omega^2|x|^2}{4}+\alp_a^2V_{\Omega}\big(\alp_a(x+y_0)
\big)\Big]\frac{\partial| R_{j,a}|^2}{\partial x_l}\\
&=\frac{1}{2}\big(\alp_a^2\mu_{j,a}+a|\tilde u_{j,a}|^2\big)\frac{\partial| R_{j,a}|^2}{\partial x_l}\quad \text{in\,\ $\R^2$}, \,\ j,l=1,2.
\end{split}
\end{equation}
Following the exponential decay \eqref{1.19}, we  calculate that  for $j=1,\,2,$
\[
\begin{aligned}
-\inte \Delta  R_{j,a}\frac{\partial R_{j,a}(x)}{\partial x_l}dx
&=-\lim_{R\to\infty}\int_{B_R(0)} \Delta  R_{j,a}\frac{\partial R_{j,a}(x)}{\partial x_l}dx\\
&=-\lim_{R\to\infty}\int_{\partial B_R(0)}\Big(\frac{\partial R_{j,a}(x)}{\partial \nu}\frac{\partial R_{j,a}(x)}{\partial x_l}- \frac {1}{2} |\nabla   R_{j,a}|^2\nu_l\Big)dS=0,
\end{aligned}
\]
where $\nu=(\nu_1,\nu_2)$ denotes the outward unit of $\partial B_{R}(0)$.
Similarly, we have
\[
\inte|\tilde u_{j,a}|^2\frac{\partial| R_{j,a}|^2}{\partial x_l}dx=\inte(R_{j,a}^2+I_{j,a}^2)\frac{\partial| R_{j,a}|^2}{\partial x_l}dx=\inte I_{j,a}^2 \frac{\partial| R_{j,a}|^2}{\partial x_l}dx,
\]
and
\[
\frac{1}{2}\inte\Big[\frac{\alp_a^4\Omega^2|x|^2}{4}+\alp_a^2V_{\Omega}(\alp_a(x+y_0)
)\Big]\frac{\partial| R_{j,a}|^2}{\partial x_l}=-\inte\frac{1}{2}\Big[\frac{\alp_a^4\Omega^2x_l}{2}+\alp_a^2\frac{\partial V_{\Omega}\big(\alp_a(x+y_0)\big)}{\partial x_l}\Big]{R}^2_{j,a},
\]
where $(x_1,x_2)=x\in \R^2$.
We then derive from \eqref{3.1} that for $l=1,\,2,$
\begin{equation*}
\arraycolsep=1.5pt
 \begin{array}{lll}
&\displaystyle\inte\frac{1}{2}\Big[\frac{\alp_a^4\Omega^2x_l}{2}+\displaystyle\alp_a^2\frac{\partial V_{\Omega}\big(\alp_a(x+y_0)\big)}{\partial x_l}\Big]{R}^2_{j,a}\\[3mm]
=&-\displaystyle\inte\alp_a^2\Omega (x^{\perp}\cdot \nabla  I_{j,a}) \frac{\partial R_{j,a}(x)}{\partial x_l}
-\displaystyle\inte\frac{a}{2}|I_{j,a}|^2\frac{\partial| R_{j,a}|^2}{\partial x_l}dx,\quad j=1,\,2.
 \end{array}\eqno{(E)_j}
\end{equation*}
Following the above identity $(E)_j$ for $j=1,\,2$, we calculate
\begin{equation}\label{3.3M}
\frac{(E)_2-(E)_1}{\|\tilde u_{2,a}(x)-\tilde u_{1,a}(x)\|_{L^{\infty}(\R^2)}},
\end{equation}
which then yields from (\ref{3.21}) and (\ref{1.24}) that
\begin{equation}\label{3.3}
\begin{split}
&\inte\frac{\alp_a^2}{2}\frac{\partial V_{\Omega}\big(\alp_a(x+y_0)\big)}{\partial x_l}({R}_{2,a}+{R}_{1,a})\eta_{1,a}\\
=&-\inte \frac{\alp_a^4\Omega^2x_l}{4} (R_{2,a}+R_{1,a}) \eta_{1,a}-\inte\Big\{\alp_a^2\Omega (x^{\perp}\cdot \nabla \eta_{2,a}) \frac{\partial R_{2,a}(x)}{\partial x_l}+\alp_a^2\Omega (x^{\perp}\cdot \nabla I_{1,a}) \frac{\partial {\eta} _{1,a}}{\partial x_l}\Big\}\\
&-\inte\frac{a}{2} {\eta} _{2,a}(  I_{1,a}+  I_{2,a})\frac{\partial | {R}_{2,a}|^2}{\partial x_l}-\inte\frac{a}{2}| {I}_{1,a}|^2\frac{\partial \big[ {\eta}_{1,a}( {R}_{1,a}+ {R}_{2,a})\big]}{\partial x_l}.
\end{split}
\end{equation}
We next prove (\ref{3.34}) by estimating all terms of \eqref{3.3}.

Following Lemma \ref{1.23}, we get that as $a\nearrow a^*$,
\begin{equation}\label{3.11}
\begin{split}
&-\inte \frac{\alp_a^4\Omega^2x_l}{4} ({R}_{2,a}+{R}_{1,a}) \eta_{1,a}dx\\
=&-\frac{\alp_a^4\Omega^2}{2\sqrt{a^*}}\inte x_lw\Big[b_0(w+x\cdot \nabla w)+\sum_{j=1}^{2}b_j\frac{\partial w}{\partial x_j}\Big]dx+o(\alp_a^4)\\
=&-\frac{\alp_a^4\Omega^2}{4\sqrt{a^*}}\inte b_lx_l\frac{\partial w^2}{\partial x_l}dx+o(\alp_a^4)\\
=&\frac{\sqrt{a^*}\Omega^2b_l\alp_a^4}{4}+o(\alp_a^4),\,\ l=1,2.
\end{split}
\end{equation}
Similarly, we derive from Proposition \ref{p-3.3} that as $a\nearrow a^*$,
\begin{equation}\label{3.4}
\begin{split}
&-\alp_a^2\Omega\inte x^{\perp}\cdot \nabla \eta_{2,a}(x) \frac{\partial R_{2,a}(x)}{\partial x_l}dx\\
=&\frac{\Omega^2\alp_a^4}{2\sqrt{a^*}}\inte x^{\perp}\cdot \nabla \big[(b_1x_2-b_2x_1)w(x)\big] \frac{\partial w(x)}{\partial x_l}dx+o(\alp_a^4)\\
=&\frac{\Omega^2\alp_a^4}{2\sqrt{a^*}}\inte \big[(b_1x_1+b_2x_2)w(x)\big] \frac{\partial w(x)}{\partial x_l}dx+o(\alp_a^4)\\
=&-\frac{\sqrt{a^*}\Omega^2b_l\alp_a^4}{4}+o(\alp_a^4),\,\ l=1,2.
\end{split}
\end{equation}
Applying Lemma \ref{prop-1} and Proposition \ref{p-3.3}, we also deduce that as $a\nearrow a^*$,
\begin{equation}\label{3.5}
\inte\alp_a^2\Omega \big(x^{\perp}\cdot \nabla  I_{1,a} \big)\frac{\partial {\eta} _{1,a}}{\partial x_l}dx=o(\alp_a^4),
\end{equation}
and
\begin{equation}\label{3.6}
\inte\frac{a}{2} {\eta} _{2,a}(  I_{1,a}+  I_{2,a})\frac{\partial | {R}_{2,a}|^2}{\partial x_l}dx=o(\alp_a^4),\quad
\inte\frac{a}{2}| {I}_{1,a}|^2\frac{\partial \big[ {\eta}_{1,a}( {R}_{1,a}+ {R}_{2,a})\big]}{\partial x_l}dx=o(\alp_a^4).
\end{equation}
Under the assumption $(V)$,  we finally conclude from (\ref{3.3})--\eqref{3.6} that as $a\nearrow a^*$,
\begin{equation}\label{3.33}
\begin{split}
 o(\alp_a^4)
&=\inte\frac{\alp_a^2}{2}\frac{\partial V_{\Omega}\big(\alp_a(x+y_0)\big)}{\partial x_l}( {R}_{2,a}+ {R}_{1,a})  \eta_{1,a}dx\\
&=\frac{\alp_a^{2+p}}{2}\inte\frac{\partial h(x+y_0)}{\partial x_l}({R}_{2,a}+{R}_{1,a})\eta_{1,a}dx+o(\alp_a^{2+p})\\
&=\frac{\alp_a^{2+p}}{\sqrt{a^*}}\inte\frac{\partial h(x+y_0)}{\partial x_l}w\Big[b_0(w+x\cdot \nabla w)+\sum_{j=1}^2b_j\frac{\partial w}{\partial x_j}\Big]dx+o(\alp_a^{2+p})\\
&=\frac{\alp_a^{2+p}}{2\sqrt{a^*}}\inte \frac{\partial h(x+y_0)}{\partial x_l}\Big[b_0(x\cdot \nabla w^2)+\sum_{j=1}^2b_j\frac{\partial w^2}{\partial x_j}\Big]dx+o(\alp_a^{2+p}),\,\ l=1,2,
\end{split}
\end{equation}
where we have used the fact that $y_0$ is the unique critical point of $H(y)$. This further implies that the claim \eqref{3.34} holds true.


\vskip 0.05truein
\noindent{\em  Step 2.}
The constant $b_0=0$ in \eqref{3d:5a}.\vskip 0.05truein

Multiplying the equation \eqref{1.3} by $ (x\cdot \nabla  R_{j,a}) $, we have for $j=1,2$,
\begin{equation}\label{3.15}
\begin{split}
& -\Delta  R_{j,a}(x\cdot \nabla  R_{j,a}) -\alp_a^2\Omega(x^{\perp}\cdot \nabla I_{j,a})(x\cdot \nabla R_{j,a})\\
=&\Big[\alp _a^2 \mu_{j,a}-\frac{\alp_a^4\Omega^2|x|^4}{4}-\alp_a^2V_{\Omega}\big(\alp_a(x+y_0)\big)\Big] R_{j,a}(x\cdot \nabla R_{j,a})\\
&+a |\tilde u_{j,a}|^2 R_{j,a}(x\cdot \nabla R_{j,a})\quad\hbox{in $\R^2$}.
\end{split}
\end{equation}
Using the integration by parts,  we note from  \eqref{1.19} that for $j=1,2,$
\begin{equation*}\arraycolsep=1.5pt\begin{array}{lll}
A_{j,a}:&=&-\displaystyle\inte \Delta  R_{j,a}(x\cdot \nabla  R_{j,a})dx\\[4mm]
&=& -\lim_{R\to\infty} \displaystyle \int_{\partial B_R(0)} \Big[\frac{\partial R_{j,a} }{\partial \nu }(x\cdot \nabla R_{j,a})-(x\cdot \nu) |\nabla \hat R_{j,a}|^2\Big]dS\\
&=&0,
\end{array}
\end{equation*}
\[
\begin{aligned}
B_{j,a}:&=\displaystyle \inte \Big[\alp _a^2\mu_{j,a}-\frac{\alp_a^4\Omega^2|x|^2}{4}-\alp_a^2V_{\Omega}\big(\alp_a(x+y_0)\big)\Big] R_{j,a}(x\cdot \nabla  R_{j,a})\\[4mm]
&=-\displaystyle\inte  R_{j,a}^2\Big[\alp _a^2\mu_{j,a}-\frac{\alp _a^4\Omega^2|x|^2}{2}-\alp_a^2V_{\Omega}\big(\alp_a(x+y_0)\big)-\frac{\alp_a^2}{2}x\cdot\nabla_x V_{\Omega}\big(\alp_a(x+y_0)\big)\Big],
\end{aligned}
\]
and
\[
\begin{aligned}
C_{j,a}:&=a \inte |\tilde u_{j,a}|^2  R_{j,a} (x\cdot \nabla  R_{j,a} ) \\[4mm]
&=\frac{a}{4}\inte (x \cdot \nabla R^4_{j,a}) +\frac{a}{2}\inte I_{j,a}^2 (x \cdot \nabla  R^2_{j,a}) \\[4mm]
&=-\frac{a}{2}\inte R^4_{j,a}
+\frac{a}{2}\inte I_{j,a}^2 (x \cdot \nabla  R^2_{j,a}) .
\end{aligned}
\]
Therefore, we derive from above that
\begin{equation}\label{3.16}
\frac{(D_{2,a}-D_{1,a})}{\|\tilde{u}_{2,a}-\tilde{u}_{1,a}\|_{L^{\infty}(\R^2)}}=\frac{(B_{2,a}-B_{1,a})+(C_{2,a}-C_{1,a})}{\|\tilde{u}_{2,a}-\tilde{u}_{1,a}\|_{L^{\infty}(\R^2)}},
\end{equation}
where $D_{j,a}$ is defined by
\[
\begin{aligned}
D_{j,a}:=-\inte\alp_a^2\Omega(x^{\perp}\cdot\nabla I_{j,a})(x\cdot \nabla R_{j,a}),\quad j=1,\,2.
\end{aligned}
\]
We next estimate the terms containing $D_{j,a}$ and  $B_{j,a}$ of (\ref{3.16}) as follows.

As for the term   containing $D_{j,a}$, applying Lemmas \ref{prop-1} and \ref{1.25}, we infer from  Proposition \ref{p-3.3} that as $a\nearrow a^*$,
\begin{equation}\label{3.20}
\begin{split}
&\quad\frac{D_{2,a}-D_{1,a}}{\|\tilde{u}_{2,a}-\tilde{u}_{1,a}\|_{L^{\infty}(\R^2)}}\\
&=-\alp_a^2\inte\Omega \big(x^{\perp}\cdot \nabla  \eta_{2,a}\big)(x\cdot \nabla R_{2,a})-\alp_a^2\inte\Omega \big(x^{\perp}\cdot \nabla  I_{1,a}\big)(x\cdot \nabla \eta_{1,a})\\
&=-\alp_a^2\inte\Omega \big(x^{\perp}\cdot \nabla  \eta_{2,a}\big)\Big[x\cdot\Big( \nabla \frac{w}{\sqrt{a^*}}\Big) +x\cdot\nabla\Big(  R_{2,a}-\frac{w}{\sqrt{a^*}}\Big)\Big]+o(\alp_a^4)\\
&=\alp_a^2\frac{\Omega}{\sqrt{a^*}}\inte [x^{\perp}\cdot \nabla \big(x\cdot \nabla w \big)  ]\eta_{2,a}+o(\alp_a^4)=o(\alp_a^4),
\end{split}
\end{equation}
where (\ref{3.21}) and (\ref{1.24}) are also used.
As for the term  containing $B_{j,a}$, we obtain from the assumption $(V)$ that as $a\nearrow a^*$,
\[
\begin{aligned}
&\quad\frac{B_{2,a}-B_{1,a}}{\|\tilde{u}_{2,a}-\tilde{u}_{1,a}\|_{L^{\infty}}}\\
&=\displaystyle  \inte\Big[\frac{\alp_a^4\Omega^2|x|^2}{2}+\alp_a^2V_{\Omega}\big(\alp_a(x+y_0)\big)+\frac{\alp_a^2}{2}(x+y_0)\cdot \nabla_x V_{\Omega}\big(\alp_a(x+y_0)\big)\Big]\\
&\quad\ \quad\quad\cdot (R_{1,a}+ R_{2,a})\eta_{1,a}dx-J_a-K_a,\\
&=\displaystyle  \inte\Big[\frac{\alp_a^4\Omega^2|x|^2}{2}+\frac{2+p}{2}\alp_a^{2+p}h(x+y_0)\Big](R_{1,a}+ R_{2,a})\eta_{1,a}dx-J_a-K_a+o(\alp_a^{2+p}),\\
\end{aligned}
\]
where the fact $x\cdot \nabla h(x)=ph(x)$ is used in the last equality. Here the terms $J_a$ and $K_a$ are defined by
\[
\begin{aligned}
J_a:=\displaystyle\frac{\alp _a^2}{2} \inte\Big[y_0\cdot \nabla_x V_{\Omega}\big(\alp_a(x+y_0)\big)\Big](R_{1,a}+ R_{2,a})\eta_{1,a}dx,
\end{aligned}
\]
and
\[
\begin{aligned}
K_a:&=\displaystyle\alp _a^2 \inte \frac{  R^2_{2,a}\mu_{2,a}-  R^2_{1,a}\mu_{1,a}}{\|\tilde{u}_{2,a}-\tilde{u}_{1,a}\|_{L^{\infty}(\R^2)}}dx.
\end{aligned}
\]
Applying \eqref{5:DDE}, one can note from the first identity of \eqref{3.33} that as $a\nearrow a^*$,
\begin{equation*}\label{s.6}
J_a=o(\alp_a^4).
\end{equation*}
Using Lemma \ref{prop-1} and Proposition \ref{p-3.3}, we derive from \eqref{4.7b} that $K_a$ satisfies
\[
\begin{aligned}
K_a:
&=\displaystyle\alp _a^2 \inte\frac{|\tilde u_{2,a}|^2\mu_{2,a}-|\tilde u_{1,a}|^2\mu_{1,a}}{\|\tilde{u}_{2,a}-\tilde{u}_{1,a}\|_{L^{\infty}(\R^2)}}
-\displaystyle\alp _a^2\inte \frac{  I^2_{2,a}\mu_{2,a}- I^2_{1,a}\mu_{1,a}}{\|\tilde{u}_{2,a}-\tilde{u}_{1,a}\|_{L^{\infty}(\R^2)}}\\
&=\frac{\alp _a^2(\mu_{2,a}-\mu_{1,a})}{\|\tilde{u}_{2,a}-\tilde{u}_{1,a}\|_{L^{\infty}(\R^2)}}-\displaystyle\alp _a^2 \inte  \eta_{2,a}( I_{1,a}+ I_{2,a})\mu_{2,a}\\
&\quad-\alp _a^2\frac{\mu_{2,a}-\mu_{1,a}} {\|\tilde{u}_{2,a}-\tilde{u}_{1,a}\|_{L^{\infty}(\R^2)}}\inte I^2_{1,a}\\
&=\frac{\alp _a^2(\mu_{2,a}-\mu_{1,a})}{\|\tilde{u}_{2,a}-\tilde{u}_{1,a}\|_{L^{\infty}(\R^2)}}+o(\alp_a^4)\ \ \mbox{as}\ \ a\nearrow a^*.
\end{aligned}
\]
It then follows from above that
\begin{equation}\label{3.18}
\begin{split}
&\quad\frac{B_{2,a}-B_{1,a}}{\|\tilde{u}_{2,a}-\tilde{u}_{1,a}\|_{L^{\infty}(\R^2)}}\\
&=\displaystyle \inte\Big[\frac{\alp _a^4\Omega^2|x|^2}{2}+\frac{2+p}{2}\alp_a^{2+p}h(x+y_0)\Big](R_{1,a}+R_{2,a})\eta_{1,a}\\
&\quad
-\frac{\alp _a^2(\mu_{2,a}-\mu_{1,a})}{\|\hat{u}_{2,a}-\hat{u}_{1,a}\|_{L^{\infty}(\R^2)}}+o(\alp_a^{2+p})\ \ \mbox{as}\ \ a\nearrow a^*.
\end{split}
\end{equation}

Applying Lemma \ref{prop-1} and Proposition \ref{p-3.3}, we obtain that as $a\nearrow a^*$,
\begin{equation}\label{3.19}
\begin{split}
&\quad\frac{C_{2,a}-C_{1,a}}{\|\tilde{u}_{2,a}-\tilde{u}_{1,a}\|_{L^{\infty}(\R^2)}}\\
&=-\frac{a}{2}\inte \frac{|\tilde u_{2,a}|^4-|\tilde u_{1,a}|^4}{\|\tilde{u}_{2,a}-\tilde{u}_{1,a}\|_{L^{\infty}(\R^2)}}+\frac{a}{2}\inte \frac{(|I_{2,a}|^4-| I_{1,a}|^4)+2(| R_{2,a}|^2|I_{2,a}|^2-| R_{1,a}|^2| I_{1,a}|^2)}{\|\tilde{u}_{2,a}-\tilde{u}_{1,a}\|_{L^{\infty}(\R^2)}}\\
&\quad
+\frac{a}{2}\inte I_{2,a}^2\big[x\cdot \nabla[ \eta_{1,a}(R_{1,a}+R_{2,a})]\big]
+\frac{a}{2}\inte  \eta_{2,a}( {I}_{1,a}+ {I}_{2,a})(x\cdot \nabla  R^2_{1,a})\\
&=-\frac{a}{2}\inte \frac{|\tilde u_{2,a}|^4-|\tilde u_{1,a}|^4}{\|\tilde{u}_{2,a}-\tilde{u}_{1,a}\|_{L^{\infty}(\R^2)}}+o(\alp^4_a).
\end{split}
\end{equation}
Note from (\ref{1.34}) that
\[
-\alp_a^2\frac{\mu_{2,a}-\mu_{1,a}}{\|\tilde{u}_{2,a}-\tilde{u}_{1,a}\|_{L^{\infty}(\R^2)}}-\frac{a}{2}\inte \frac{|\tilde u_{2,a}|^4-|\tilde u_{1,a}|^4}{\|\tilde{u}_{2,a}-\tilde{u}_{1,a}\|_{L^{\infty}(\R^2)}}dx\equiv 0.
\]
Following this identity, we then conclude from \eqref{3.16}--\eqref{3.19} that as $a\nearrow a^*$,
\begin{equation}\label{4.28a}
\begin{aligned}
o(\alp_a^4)=&
\frac{(B_{2,a}-B_{1,a})+(C_{2,a}-C_{1,a})}{\|\tilde{u}_{2,a}-\tilde{u}_{1,a}\|_{L^{\infty}(\R^2)}}\\
=&\displaystyle  \inte\Big[\frac{\Omega^2\alp _a^4|x|^2}{2}+\frac{(2+p)\alp_a^{2+p}}{2}h(x+y_0)\Big]( R_{1,a}+ R_{2,a}) \eta_{1,a}dx+o(\alp_a^{2+p})\\
=&\displaystyle\frac{(2+p)\alp _a^{2+p}}{\sqrt{a^*}} \inte h(x+y_0)w \eta_{1}dx+\displaystyle\frac{\Omega^2\alp _a^{4}}{\sqrt{a^*}} \inte |x|^2w \eta_{1}dx+o(\alp_a^{2+p}).
\end{aligned}
\end{equation}
Since $y_0$ is the unique critical point of $H(y)$, we finally derive from (\ref{4.28a}) that if $1<p<2$,
\[
\begin{aligned}
0&=\inte h(x+y_0)w \eta_{1}dx\\
 &=b_0\inte h(x+y_0)w\big(w+x\cdot \nabla w\big)dx+\sum_{i=1}^2b_i\inte \frac 12h(x+y_0)\frac{\partial_iw^2}{\partial x_i}dx\\
 &=b_0\inte h(x+y_0)w^2dx-\frac{b_0}{2}\inte \Big(x\cdot \nabla h(x+y_0)w^2+2h(x+y_0)w^2\Big)dx\\
 &=-\frac{b_0}{2}\inte \Big(\big(x+y_0\big)\cdot \nabla h(x+y_0)w^2\Big)dx\\
 &=-\frac{pb_0}{2}\inte h(x+y_0)w^2dx,
\end{aligned}
\]
which further implies that $b_0=0$. Similarly, the above conclusion also holds true in the case $p=2$, and the claim $b_0=0$ is therefore proved.

\vskip 0.05truein
\noindent{\em  Step 3.} The constants $b_1=b_2=0$.\vskip 0.05truein

By Step 2, we deduce from \eqref{3.34} that
\begin{equation}\label{A.35}
\sum_{j=1}^2b_j\inte \frac{\partial h(x+y_0)}{\partial x_l}\frac{\partial w^2}{\partial x_j}dx=0,\,\ l=1,2.
\end{equation}
Using the non-degeneracy assumption of the critical point $y_0$, we then derive from (\ref{A.35}) that $b_1=b_2=0$, and the proof of Step 3 is therefore complete.

\vskip 0.05truein

By the exponential decay of Lemma \ref{1.23}, we  obtain that $\eta_a\to \eta_0=\eta_1+i\eta_2\not\equiv 0$ uniformly in $C^1(\mathbb{R}^2)$ as $a \nearrow a^*$, due to the fact that $\|\eta _a\|_{L^\infty}\equiv 1$. However, we conclude from Step 2 and Step 3 that $\eta_0 \equiv 0$, a contradiction. This completes the proof of Theorem \ref{thm1.2}. \qed


 \vskip 0.16truein
\noindent {\bf Acknowledgements:}
The authors are  very grateful to the referee for many valuable suggestions which lead to the great improvements of the present paper.

\end{document}